\documentclass[review,onefignum,onetabnum]{siamart220329}
\usepackage{braket,amsfonts,amsopn}
\usepackage{amsmath, amssymb, amsfonts}
\usepackage{graphicx,epstopdf} 
\usepackage[caption=false]{subfig} 
\usepackage{enumitem}
\usepackage{color,xcolor}
\graphicspath{{figures/}}

\usepackage{lipsum}
\usepackage{amsfonts}
\usepackage{graphicx}
\usepackage{epstopdf}
\usepackage{algorithmic}
\ifpdf
  \DeclareGraphicsExtensions{.eps,.pdf,.png,.jpg}
\else
  \DeclareGraphicsExtensions{.eps}
\fi

\newsiamremark{remark}{Remark}
\newsiamremark{hypothesis}{Hypothesis}
\crefname{hypothesis}{Hypothesis}{Hypotheses}
\newsiamthm{claim}{Claim}

\headers{An Example Article}{D. Doe, P. T. Frank, and J. E. Smith}

\title{An Example Article\thanks{Submitted to the editors DATE.
\funding{This work was funded by the Fog Research Institute under contract no.~FRI-454.}}}

\author{Dianne Doe\thanks{Imagination Corp., Chicago, IL 
  (\email{ddoe@imag.com}, \url{http://www.imag.com/\string~ddoe/}).}
\and Paul T. Frank\thanks{Department of Applied Mathematics, Fictional University, Boise, ID 
  (\email{ptfrank@fictional.edu}, \email{jesmith@fictional.edu}).}
\and Jane E. Smith\footnotemark[3]}

\usepackage{amsopn}

\ifpdf
\hypersetup{
  pdftitle={An Example Article},
  pdfauthor={D. Doe, P. T. Frank, and J. E. Smith}
}
\fi

\usepackage{lineno}          
\nolinenumbers              

\definecolor{darkblue}{rgb}{0,0,0.5}
\definecolor{darkgreen}{rgb}{0,0.6,0}
\definecolor{mymagenta}{RGB}{238,0,238}

\begin{document}

\title{Algebraic Riccati Tensor Equations with Applications in Multilinear Control Systems\thanks{Submitted to the editors
			DATE. ??? \funding{This work was supported in part by Innovation Program for Quantum Science and Technology 2023ZD0300600, National Natural Science Foundation of China under Grant Nos. U24A2001, 12271108,  62173288 and 62473117, the Science and Technology Commission of Shanghai	Municipality under Grant 23JC1400501,  Guangdong Provincial Quantum Science Strategic Initiative (No. GDZX2200001), Hong Kong Research Grant Council (RGC) under Grant No. 15213924, and the CAS AMSS-PolyU Joint Laboratory of Applied Mathematics}}}
\author{Yuchao Wang\thanks{School of Mathematical Sciences, Fudan University, Shanghai, P. R. of China. (\email{yuchaowang21@m.fudan.edu.cn.})} \and Yimin Wei\thanks{School of Mathematical Sciences and  Key Laboratory of Mathematics for Nonlinear Sciences, Fudan University, Shanghai, P. R. of China. (\email{ymwei@fudan.edu.cn.})} \and Guofeng Zhang\thanks{Corresponding author (G. Zhang). Department of Applied Mathematics, The Hong Kong Polytechnic University, Kowloon 999077, Hong Kong, China, and Research Institute for Quantum Technology, The Hong Kong Polytechnic University,  Hong Kong, China, and The Hong Kong Polytechnic University Shenzhen Research Institute, Shenzhen, 518057, China.  (\email{guofeng.zhang@polyu.edu.hk.})}
\and Shih Yu Chang\thanks{Department of Applied Data Science, San Jose State University, San Jose, CA, USA. (\email{shihyu.chang@sjsu.edu})}}

\headers{Algebraic Riccati Tensor Equations and Multilinear Control}{Yuchao Wang, Yimin Wei, Guofeng Zhang and Shih Yu Chang}

\maketitle

\begin{abstract} In a recent paper by Chen et al. \cite{chen2021multilinear}, the authors initiated the control-theoretic study of a class of discrete-time multilinear time-invariant (MLTI) control systems, where system states, inputs, and outputs are all tensors endowed with the Einstein product. They established criteria for fundamental system-theoretic notions such as stability, reachability, and observability through tensor decomposition.
Building on this new research direction, the purpose of our paper is to extend the study to continuous-time MLTI control systems. Specifically, we define Hamiltonian tensors and symplectic tensors, and we establish the Schur-Hamiltonian tensor decomposition and the symplectic tensor singular value decomposition (SVD). Based on these concepts, we propose the algebraic Riccati tensor equation (ARTE) and demonstrate that it has a unique positive semidefinite solution if the system is stabilizable and detectable.
To find numerical solutions to the ARTE, we introduce a tensor-based Newton method. Additionally, we establish the tensor versions of the bounded real lemma and the small gain theorem. A first-order robustness analysis of the ARTE is also conducted. Finally, we provide a numerical example to illustrate the proposed theory and algorithms.
\end{abstract}



\section{Introduction}

As natural extension of vectors (1D arrays) and matrices (2D arrays) to higher dimensional arrays,  tensors are flexible in representing data and able to characterize higher-order interactions. Thus, tensors have found wide applications in a multitude of fields such as signal processing, social and biological networks, scientific computing, and quantum physics \cite{kolda2009tensor,CZRA09,LPV13,QCC18,orus19,DCL+22,ZD22,Nie23}. For example, in biological networks, the classic Lotka-Volterra model has recently been extended to including higher-order interactions among species in order to better describe nonlinear interactions, therefore, the system evolution is governed by  tensors \cite{LS19}.  In \cite{CSJC23}, stability analysis of the higher-order Lotka-Volterra model has been conducted. 

The Einstein product for tensors is widely used in various fields such as astrophysics \cite{einstein1916foundation}, elastodynamics \cite{itskov2000theory}, data processing \cite{chang2021multi}, and hypergraphs \cite{pickard2023kronecker,CZJC24}. Due to their important applications in numerical solvers for partial differential equations, continuum physics, and engineering problems, tensor equations via the  Einstein product have attracted a lot of attention. For example, tensor Krylov subspace methods were designed for solving the Sylvester tensor equation \cite{ali2016krylov,huang2021numerical, hajarian2020conjugate}.

In \cite{surana2016dynamic, rogers2013multilinear, ding2018tensor} the authors studied a class of discrete-time MLTI  systems of the form
\begin{align}  \label{MLTI-mode}
		\mathcal{X}_{t+1} =\mathcal{X}_{t} \times_1 A_1 \times_2 \cdots \times_N A_N,  \ \ 
		\mathcal{Y}_t = \mathcal{X}_t \times_1 C_1 \times_2\cdots \times _N C_N,
\end{align}
where $\mathcal{X}_t$ is the state tensor and  $\mathcal{Y}_t$ is the output tensor at time instant $t$, $A_n$ and $C_n$ are constant matrices, and  ``$\times_n$'' is the mode product commonly used in tensor operations (\cref{def:mode product}). These authors numerically demonstrated the superior performance of the MLTI system \cref{MLTI-mode}  compared to the classical discrete-time linear time-invariant (LTI) systems in skeleton-based human behavior recognition \cite{ding2018tensor}, the modeling of regional climate change \cite{rogers2013multilinear} and dynamic texture videos \cite{surana2016dynamic}.  Multipartite quantum states can be nicely represented by tensors \cite[Chapter 9]{QCC18}. For example, it is shown in \cite{Z14, Z17, ZD22} that the steady-state response of a stable quantum linear passive system to a multi-channel multi-photon input state can be compactly represented in terms of tensors; essentially speaking, the output tensor is the higher-order SVD (HOSVD) of the input tensor (served as the core tensor) and the impulse response matrix function (or equivalently the transfer matrix function) of the system. If the input tensor passes through a series of quantum systems sequentially, the dynamics are described by the MLTI system \cref{MLTI-mode} with time-varying matrices $A_n$.

In recent works \cite{chen2019multilinear, chen2021multilinear}, Chen et al. generalized the MLTI system \cref{MLTI-mode} using the Einstein product and even-order paired tensors, and proposed the following discrete-time MLTI control systems
\begin{align} \label{MLTI-Einstein}
		\mathcal{X}_{t+1} = \mathcal{A}* \mathcal{X}_{t} + \mathcal{B}*\mathcal{U}_t, \ \ 
		\mathcal{Y}_t = \mathcal{C}*\mathcal{X}_t,
\end{align}
in which ``$*$'' stands for the Einstein product (\cref{def:einstein}),  the input $\mathcal{U}_t$, output $\mathcal{Y}_t$, and state $\mathcal{X}_t$ are discrete-time $N$th-order tensorial processes, while $\mathcal{A}, \mathcal{B}, \mathcal{C}$ are  constant $2N$th-order paired tensors. The system \cref{MLTI-mode}  can be obtained from the  system \cref{MLTI-Einstein}  when $\mathcal{A} = A_1\circ A_2 \circ \cdots \circ A_N, \mathcal{C} = C_1\circ C_2\circ \cdots \circ C_N$  and $\mathcal{B}=\mathcal{O}$ (the zero tensor).  Stability, reachability and observability for discrete-time MLTI systems \cref{MLTI-Einstein} are established under the framework of tensor algebraic theory, which are generalizations of the  control theory of discrete-time LTI systems.

However, in many scenarios, tensorial time series such as those in computer vision and financial data analysis are tensor-valued processes on continuous time intervals, and the observed multidimensional data are discretizations of the underlying continuous-time tensor-valued process \cite{bai2021tensorial}. The tensor-valued differential equations \cite{bai2022tensor, bai2021tensorial} are proposed for modeling the underlying continuous-time tensorial processes in multimodal dynamical systems, which successfully preserve the data spatial and latent continuous sequential information. 

Inspired by the above studies, in this paper we study a class of continuous-time MLTI systems of the form
\begin{align} \label{MLTI-continuous}
		\dot{\mathcal{X}}(t) = \mathcal{A}* \mathcal{X}(t) + \mathcal{B}*\mathcal{U}(t), \ \ 
		\mathcal{Y}(t) = \mathcal{C}*\mathcal{X}(t)+ \mathcal{D}*\mathcal{U}(t),
\end{align}
where the state $\mathcal{X}(t)$, input $\mathcal{U}(t)$, and output $\mathcal{Y}(t) $ are continuous-time  $N$th-order tensor processes,  while $\mathcal{A},\mathcal{B}, \mathcal{C}, \mathcal{D}$ are constant $2N$th-order paired tensors. Clearly, system \cref{MLTI-continuous} includes the system $\dot{\mathcal{X}}_t =\mathcal{X}_{t} \times_1 A_1 \times_2 \cdots \times_N A_N,  \ \mathcal{Y}_t = \mathcal{X}_t \times_1 C_1 \times_2\cdots \times _N C_N$ as a special case, which is  the continuous-time counterpart of system \cref{MLTI-mode}.

In particular, we study the following algebraic Riccati tensor equation (ARTE)
\begin{equation} \label{ARTE}
    \mathcal{A}^{\mathrm{H}}* \mathcal{E} + \mathcal{E}* \mathcal{A} - \mathcal{E}*\mathcal{G}* \mathcal{E} + \mathcal{K} = \mathcal{O},
 \end{equation}
where $\mathcal{A, G, K}$ and the unknown tensor $\mathcal{E}$ are even-order paired tensors. When the order $N=1$, the ARTE \cref{ARTE} reduces to the well-known algebraic Riccati equation (ARE) 
\begin{equation} \label{ARE}
    A^{\rm H}E + EA - EGE + K = O.
\end{equation}
The linear-quadratic optimal control problem is closely related to the positive semidefinite solution of the ARE \cref{ARE} \cite{lancaster1995algebraic, ZDG96}. Analogously, the  minimizer of the multilinear quadratic optimal control problem $\min_{\mathcal{U}} \int_{0}^{\infty} \left( \mathcal{X}^{\rm H}(t)*\mathcal{G}*\mathcal{X}(t)+  \mathcal{U}^{\rm H}(t) *\mathcal{K}* \mathcal{U}(t) \right)dt$ for the continuous-time MLTI system \cref{MLTI-continuous} is  
$\mathcal{U}(t) = -\mathcal{B}^{\rm H}*\mathcal{E}*\mathcal{X}(t)$, where $\mathcal{E}$ is  the positive semidefinite solution of the ARTE  \cref{ARTE} 
with $\mathcal{G} = \mathcal{B} * \mathcal{B}^{\mathrm{H}}$ and $\mathcal{K} = \mathcal{C}^{\mathrm{H}} * \mathcal{C}$. 

The contribution of this paper is two-fold. On one hand, we continue to develop theories of  even-order paired tensors  endowed with the Einstein product. Specifically, we define symplectic tensors (\cref{def:sym}) and Hamiltonian tensors (\cref{def:ham}), and propose 
the Schur-Hamiltonian tensor decomposition (\cref{thm:Schur_Hamiltonian}) and the symplectic tensor SVD (\cref{Thm: SymSVD}). We define the tensor Vec operation (\cref{Def: tensorVec}) and discuss its properties (\cref{prop: TVO1,prop: TVO2} and \cref{Lem: KronUeigen}) with the tensor Kronecker product (\cref{Pro: TensorKron} and 
\cref{Pro: tensorKronQR}). Moreover, we discuss fast computational methods of these operations for some structured tensors (\cref{Pro: 3.12,Pro: KronSum,Pro: rankoneSchur,Pro: rankonenorm,prop:kron}). We also present and study a Lyapunov tensor equation (\cref{lem: LTE}). On the other hand, we develop continuous-time MLTI control theory. In particular, we investigate the existence and uniqueness of the positive semidefinite solution to the ARTE \cref{ARTE}, see  \cref{Thm: ARTEunique}, and conduct a first-order robustness analysis (\Cref{subsec:pertub}). Moreover, we present a tensor version of the bounded real lemma (\cref{Pro: H_norm,thm:BRL}) and the small gain theorem (\cref{lem:small gain} and \cref{thm:small gain}). 

The paper is organized as follows. In \Cref{sec: pre}, we recall related tensor notations and results. In \Cref{sec. 3}, we introduce theories of Hamiltonian and  symplectic tensors, define the tensor Vec operation for even-order paired tensors and discuss its operational properties with tensor Kronecker product, and study the Sylvester tensor equation and Lyapunov tensor equation. In \Cref{sec. 4}, the basic control theory for continuous-time MLTI systems is introduced, the unique existence of the positive semidefinite solution of the ARTE \cref{ARTE} is established, and its first-order perturbation analysis is conducted. The tensor version of the bounded real lemma and small gain theorem is proven. \Cref{sec. 5} derives the Newton method for the  ARTE \cref{ARTE}; especially for the ARTE with coefficient tensors of generalized CPD format, the computational cost can be saved.  Moreover, a numerical example is used to demonstrate the  proposed theory and algorithms. \Cref{sec:con} gives the concluding remarks.

\section{Preliminaries}
\label{sec: pre}

	{\it Notation.}   $\iota=\sqrt{-1}$. $\mathbb{C}^{I_1\times \cdots \times I_N}$ denotes the set of all $N$th-order tensors of dimensions $(I_1, \dots,I_N)$ over the complex number field $\mathbb{C}$, where $I_n$ is the dimension of the $n$th mode of the tensor. We use calligraphic capital letters such as $\mathcal{X}$,  capital letters such as $X$, and boldface lowercase letters such as $\mathbf{x}$ to denote tensors, matrices, and vectors, respectively. Given two vectors of indices  $\mathbf{i}=\left( i_1,\dots,i_N \right)$ and  $\mathbf{j} = ( j_1,\dots,j_N )$, let  $\mathbf{i}\wedge \mathbf{j} = ( i_1,j_1,\dots, i_N,j_N )$ . The colon notation `:'  is used to specify index range (as in MATLAB). For example, $\mathcal{X}_{i_1\cdots i_N}$ is the entry of the tensor $\mathcal{X}$  specified by the index vector  $\mathbf{i}$, while $\mathcal{X}_{:i_2\cdots i_N}$ is the vector obtained by fixing indices other than the first one. Given two vectors of mode  dimensions $\mathbf{I} = ( I_1,\dots,I_N )$ and  $\mathbf{J} = ( J_1, \dots,J_N )$, their  Hadamard product $\mathbf{I}\odot \mathbf{J} = ( I_1J_1,\dots, I_NJ_N )$, namely the elementwise product of $\mathbf{I}$ and $\mathbf{J}$.  Denote $\left | \mathbf{I} \right |=  {\textstyle \prod_{n=1}^{N}}I_n $.  Finally, the vector-valued inequality $\mathbf{1} \le \mathbf{i} \le \mathbf{I}$ is understood as $1\le i_n \le I_n$ for all $n = 1,\dots, N$. 
	
\subsection{Even-order paired tensors}
Even-order paired tensors were originally proposed by Huang and Qi \cite{huang2018positive},  based on which Chen et al. \cite{chen2019multilinear, chen2021multilinear}  developed  the control theory for the discrete-time  MLTI systems \cref{MLTI-Einstein} via the Einstein product. 
\begin{definition}[Even-order paired tensors \cite{chen2019multilinear,chen2021multilinear,huang2018positive}] \label{Def: EoP_tensor}
	Even-order paired tensors are $2$Nth-order tensors with entries specified by a pairwise index notation, i.e., $\mathcal{A}_{i_1j_1\cdots i_Nj_N}$ with paired indices $(i_n, j_n)$ for each $n= 1,2,\dots,N$. The $(2n-1)$th mode of $\mathcal{A}$ is called the $n$-mode row, while the $2n$th mode of $\mathcal{A}$ is called the $n$-mode column.
\end{definition}

\begin{definition}[Einstein product \cite{einstein1916foundation,Brazell2013,chen2021multilinear}] \label{def:einstein}
	Given two even-order paired tensors $\mathcal{A}\in \mathbb{C}^{I_1\times J_1 \times \cdots \times I_N\times J_N}$ and $\mathcal{B}\in \mathbb{C}^{J_1\times K_1\times \cdots \times J_N\times K_N}$, the Einstein product $\mathcal{A}*\mathcal{B}$ generates a tensor in $ \mathbb{C}^{I_1\times K_1\times  \cdots \times I_N\times K_N}$,  whose entries are
	\[ \left( \mathcal{A}*\mathcal{B} \right)_{i_1k_1\cdots i_Nk_N} = \sum_{j_1=1}^{J_1} \cdots \sum_{j_N=1}^{J_N}\mathcal{A}_{i_1j_1\cdots i_Nj_N}\mathcal{B}_{j_1k_1\cdots j_Nk_N}. \]
\end{definition}

\begin{remark}
\label{Remark: 2.1}
    We regard an $N$th-order tensor $\mathcal{Y}$ as a $2N$th-order paired tensor whose dimensions of all mode columns are 1. As a result, the definition of the Einstein product for even-order paired tensors can be adopted. For example, if $\mathcal{Y} \in \mathbb{C}^{J_1 \times \cdots \times J_N}$, the Einstein product $\mathcal{A}*\mathcal{Y} \in \mathbb{C}^{I_1 \times \cdots \times I_N}$ is  $\left( \mathcal{A}*\mathcal{Y} \right)_{i_1\cdots i_N} = \sum_{j_1=1}^{J_1} \cdots \sum_{j_N=1}^{J_N}\mathcal{A}_{i_1j_1\cdots i_Nj_N}\mathcal{Y}_{j_1\cdots j_N}$.	If $\mathcal{Y} \in \mathbb{C}^{I_1 \times \cdots \times I_N}$, the Einstein product $\mathcal{Y}^{\rm H} *\mathcal{A} \in \mathbb{C}^{J_1 \times \cdots \times J_N}$ is $	\left( \mathcal{Y}^{\rm H}*\mathcal{A} \right)_{j_1\cdots j_N} = \sum_{i_1=1}^{I_1} \cdots \sum_{i_N=1}^{I_N} \overline{\mathcal{Y}_{i_1\cdots i_N}}\mathcal{A}_{i_1j_1\cdots i_Nj_N}$, where $\mathcal{Y}^{\rm H}$ denotes the complex conjugate transpose of $\mathcal{Y}$ with entries $(\mathcal{Y}^{\rm H})_{i_1\cdots i_N} = \overline{\mathcal{Y}_{i_1\cdots i_N}}$, and we regard $\mathcal{Y}^{\rm H}$ as a $2N$th-order paired tensor while all its dimensions of mode rows are 1.
\end{remark}

Tensor unfolding, the process of rearranging the entries of a tensor into a matrix or vector, is a frequently used operation in tensor analysis and numerical computations \cite{Brazell2013, kolda2009tensor}. Given a sequence of tensor dimensions $\mathbf{I} = \left( I_1,\dots, I_N \right)$ and a corresponding vector of indices $\mathbf{i} = \left(i_1,\dots, i_N\right)$, define 
\begin{equation} \label{eq:ivec}
  \mathrm{ivec}(\mathbf{i}, \mathbf{I}) := i_1+\sum_{k=2}^{N} (i_k-1)\prod_{j=1}^{k-1}I_j.  
\end{equation}

\begin{definition}[Tensor unfolding \cite{Brazell2013,chen2021multilinear}]\label{def:phi}
    The unfolding process of an even-order paired tensor $\mathcal{A}\in \mathbb{C}^{I_1\times J_1 \times \cdots \times I_N\times J_N}$ is the isomorphic map $\phi : \mathbb{C}^{I_1\times J_1 \times \cdots \times I_N\times J_N} \to \mathbb{C}^{|\mathbf{I}|\times |\mathbf{J}|}$ with $\phi(\mathcal{A}) = A$ defined componentwise as $ \mathcal{A}_{i_1j_1\cdots i_Nj_N} \overset{\phi}{\rightarrow} A_{\mathrm{ivec}(\mathbf{i}, \mathbf{I}) \mathrm{ivec}(\mathbf{j}, \mathbf{J})}$.
\end{definition}

Based on tensor unfolding, some notations for even-order paired tensors analogous to matrices are introduced in \cite{chen2021multilinear}, which are listed below.

\begin{itemize}[itemsep = 0pt, topsep=0pt]
	\item An even-order paired tensor $\mathcal{D}\in \mathbb{C}^{I_1\times I_1 \times \cdots \times I_N \times I_N}$ is {\it diagonal} if all its entries are $0$ except $\mathcal{D}_{i_1i_1\cdots i_Ni_N}$. Moreover, if all the diagonal entries $\mathcal{D}_{i_1i_1\cdots i_Ni_N} = 1$, then $\mathcal{D}$ is the {\it identity tensor} denoted by $\mathcal{I}_{\mathbf{I}}$, simply written as $\mathcal{I}$ when its dimensions are clear in the context.  Clearly, $I=\phi(\mathcal{I})$ is an identity matrix.
	\item Given $\mathcal{A}\in \mathbb{C}^{I_1\times I_1\times \cdots \times I_N\times I_N}$, if there exists a tensor $\mathcal{B}$ of the same size such that $\mathcal{A}*\mathcal{B}=\mathcal{B}*\mathcal{A}=\mathcal{I}$, then $\mathcal{B}$ is called the {\it inverse} of $\mathcal{A}$ denoted by $\mathcal{A}^{-1}$. An even-order paired tensor $\mathcal{U}\in \mathbb{C}^{I_1\times I_1\times \cdots \times I_N \times I_N}$ is {\it unitary} if $\mathcal{U}^{\mathrm{H}}*\mathcal{U}=\mathcal{U}*\mathcal{U}^{\mathrm{H}}=\mathcal{I}$.
	\item For an even-order paired tensor $\mathcal{A} \in \mathbb{C}^{I_1\times J_1 \times \cdots \times I_N\times J_N}$, the tensor $\mathcal{B}\in \mathbb{C}^{J_1\times I_1\times \cdots J_N\times I_N}$ is called the {\it transpose} of $\mathcal{A}$ denoted by $\mathcal{A}^{\top}$, if $\mathcal{B}_{j_1i_1\cdots j_Ni_N} = \mathcal{A}_{i_1j_1\cdots i_Nj_N}$, and is called the {\it complex conjugate transpose} of $\mathcal{A}$ denoted by $\mathcal{A}^{\mathrm{H}}$, if $\mathcal{B}_{j_1i_1\cdots j_Ni_N}=\overline{\mathcal{A}_{i_1j_1\cdots i_Nj_N}}$. $\mathcal{A}$ is said to be {\it Hermitian} if $\mathcal{A}=\mathcal{A}^{\mathrm{H}}$.
	\item The {\it unfolding rank} of an even-order paired tensor $\mathcal{A}\in \mathbb{C}^{I_1\times I_1 \times \cdots \times I_N\times I_N}$ is defined as $\mathrm{rank}_U(\mathcal{A}) := \mathrm{rank}(\phi(\mathcal{A}))$, and the {\it unfolding determinant} of $\mathcal{A}$ is defined as $\mathrm{det}_U(\mathcal{A}) := \mathrm{det}(\phi(\mathcal{A}))$.
	\item  An even-order paired Hermitian tensor $\mathcal{A}\in \mathbb{C}^{I_1\times I_1\times \cdots \times I_N\times I_N}$ is {\it positive semidefinite} if  $\mathcal{Y}^{\mathrm{H}}*\mathcal{A}*\mathcal{Y} = \sum_{\mathbf{1}\le \mathbf{i}\le \mathbf{I}}\sum_{\mathbf{1}\le \mathbf{j}\le \mathbf{I}} \overline{\mathcal{Y}_{i_1\cdots i_N}} \mathcal{A}_{i_1j_1 \cdots i_Nj_N} \mathcal{Y}_{j_1\cdots j_N} \ge 0$ for any nonzero $N$th-order tensor $\mathcal{Y} \in \mathbb{C}^{I_1\times \cdots \times I_N}$, and {\it positive definite} if the inequality is strict.
\end{itemize}
Finally, we say that an even-order paired tensor $\mathcal{T}\in \mathbb{C}^{I_1\times I_1\times \cdots \times I_N\times I_N}$ is {\it upper triangular} if $\mathcal{T}_{i_1j_1\cdots i_Nj_N} = 0$ for all $\mathrm{ivec}(\mathbf{i}, \mathbf{I}) > \mathrm{ivec}(\mathbf{j}, \mathbf{I})$.

\subsection{Tensor concatenation and blocks}\label{sec:tensor_blocking}
Like block matrices,  block tensors are tensors whose entries are themselves tensors. Here we adopt a compact concatenation approach \cite{chen2019multilinear, chen2021multilinear}  to construct block tensors.  
\begin{definition}[$n$-mode block tensor \cite{chen2021multilinear}]
\label{Def: n-modeblocktensor}
	Let $\mathcal{A}, \mathcal{B} \in \mathbb{C}^{I_1\times J_1 \times \cdots \times I_N \times J_N}$. For each $n = 1,\dots, N$, the $n$-mode row block tensor  concatenated by $\mathcal{A}$ and $\mathcal{B}$, denoted by $\left[\begin{smallmatrix}
		\mathcal{A} & \mathcal{B}
	\end{smallmatrix}\right]_n \in \mathbb{C}^{I_1\times J_1 \times \cdots \times I_n \times 2J_n \times \cdots \times I_N \times J_N}$, is defined element-wisely  as
	\begin{equation*}
		(\begin{bmatrix}
			\mathcal{A} & \mathcal{B}
		\end{bmatrix}_n)_{i_1j_1\cdots i_Nj_N} = \left\{
		\begin{aligned}
			&\mathcal{A}_{i_1j_1\cdots i_nj_n\cdots i_Nj_N},\ \quad \quad i_k = 1,\dots , I_k, j_k = 1,\dots , J_k,\forall k, \\
			&\mathcal{B}_{i_1j_1\cdots i_n(j_n -J_n)\cdots i_Nj_N}, \ i_k = 1,\dots , I_k,\forall k, j_k = 1,\dots, J_k \\
			&\quad \quad \quad \ \ \quad \quad \quad \quad \quad \quad \mathrm{for}\ k\ne n \ \mathrm{and}\ j_n = J_n + 1, \dots, 2J_n.
		\end{aligned}
	\right.
	\end{equation*}
	The $n$-mode column block tensor  is $ \left[\begin{smallmatrix}
		\mathcal{A} \\
		\mathcal{B}
	\end{smallmatrix}\right]_n := \left[\begin{smallmatrix}
		\mathcal{A}^{\top} & \mathcal{B}^{\top}
	\end{smallmatrix}\right]_n^{\top}$.
\end{definition}

We also denote by 
$\left[\begin{smallmatrix}
	\mathcal{A} &  \mathcal{B} \\
	\mathcal{C} &  \mathcal{D}
\end{smallmatrix}\right]_n =\left[\begin{smallmatrix}
	\left[\begin{smallmatrix}
		\mathcal{A} & \mathcal{B}
\end{smallmatrix}\right]_n\\
	\left[\begin{smallmatrix}
		\mathcal{C} & \mathcal{D}
	\end{smallmatrix}\right]_n
\end{smallmatrix}\right]_n$, the $n$-mode block tensor concatenated by the $n$-mode row block tensors $\left[\begin{smallmatrix}
	\mathcal{A} & \mathcal{B}
\end{smallmatrix}\right]_n$ and $\left[\begin{smallmatrix}
	\mathcal{C} & \mathcal{D}
\end{smallmatrix}\right]_n$.
\begin{remark}
	Given two $N$th-order tensors $\mathcal{X} \in \mathbb{C}^{J_1\times \cdots \times J_N}$, $\mathcal{Y} \in \mathbb{C}^{J_1 \times \cdots \times J_N}$, the $n$-mode column block tensor $\left[\begin{smallmatrix}
		\mathcal{X} \\
		\mathcal{Y}
	\end{smallmatrix}\right]_n \in \mathbb{C}^{J_1 \times \cdots \times J_{n-1}\times 2J_n\times J_{n+1}\times \cdots \times J_N}$ can also be defined as in \cref{Def: n-modeblocktensor}, by regarding $\mathcal{X}, \mathcal{Y}$ as even-order paired tensors (\cref{Remark: 2.1}).  
\end{remark}

    Under the Einstein product, block tensors enjoy  properties similar to their matrix counterparts.
\begin{proposition}[\cite{chen2021multilinear}] \label{Pro: blocktensor}
	Let $\mathcal{A}, \mathcal{B} \in \mathbb{C}^{I_1\times J_1 \times \cdots \times I_N \times J_N}$, $ \mathcal{C}, \mathcal{D} \in \mathbb{C}^{J_1\times I_1 \times \cdots \times J_N\times I_N}$. The following properties of block tensors hold for all $n=1,\dots,N$.
	\begin{enumerate} \rm
		\item
		$	\begin{bmatrix}
			\mathcal{P}* \mathcal{A} & \mathcal{P}*\mathcal{B}
		\end{bmatrix}_n = \mathcal{P}* \begin{bmatrix}
			\mathcal{A} & \mathcal{B}	\end{bmatrix}_n,\quad for \  \mathcal{P} \in \mathbb{C}^{L_1\times I_1 \times \cdots \times L_N \times I_N}$. 
		\item 
		$	\begin{bmatrix}
			\mathcal{C}* \mathcal{Q} \\
			\mathcal{D}* \mathcal{Q}
		\end{bmatrix}_n = \begin{bmatrix}
			\mathcal{C} \\
			\mathcal{D}
		\end{bmatrix}_n* \mathcal{Q}, \quad   for \  \mathcal{Q} \in \mathbb{C}^{I_1\times L_1 \times \cdots \times I_N \times L_N}$.
		\item
		$		\begin{bmatrix}
			\mathcal{A} &  \mathcal{B} 
		\end{bmatrix}_n* \begin{bmatrix}
			\mathcal{C} \\
			\mathcal{D}
		\end{bmatrix}_n = 
		\mathcal{A}* \mathcal{C} + \mathcal{B}* \mathcal{D}. $
	\end{enumerate}
\end{proposition}

Given even-order paired tensors $\mathcal{X}_{k} \in \mathbb{C}^{I_1 \times J_1 \times \cdots \times I_N\times J_N}$ $(k= 1,\ldots, K_n)$, one can apply \cref{Def: n-modeblocktensor} successively to create a bigger even-order $n$-mode row block tensor $\begin{bmatrix}
	\mathcal{X}_1 & \cdots & \mathcal{X}_{K_n}
\end{bmatrix}_n \in \mathbb{C}^{I_1\times J_ 1\times \cdots \times I_n \times J_nK_n \times  \cdots \times I_N\times J_N}$. Based on these constructions,  we give the definitions of  mode block tensors and tensor blockings, which will be used in \Cref{Section: 3.2}. 

\begin{definition}[Mode block tensor \cite{chen2019multilinear,chen2021multilinear}] 
	\label{Def: modeblocktensor}
	Given $|\mathbf{K}| = K_1K_2\cdots K_N$ even-order paired tensors $\mathcal{X}_m \in \mathbb{C}^{I_1 \times J_1 \times \cdots \times I_N\times J_N}$ $(m=1,\dots, |\mathbf{K}|)$, the mode row block tensor $\mathcal{X}$ of size $I_1 \times J_1K_1 \times \cdots \times I_N \times J_NK_N$ can be constructed by tensor concatenations of \cref{Def: n-modeblocktensor} sequentially:
	\begin{enumerate}
		\item Divide the $|\mathbf{K}|$ tensors $\{ \mathcal{X}_{1}, \mathcal{X}_2, \dots, \mathcal{X}_{|\mathbf{K}|} \}$ sequentially into $K_2\cdots K_N$ groups $\{ \mathcal{X}_1, \dots, \mathcal{X}_{K_1} \}$, $\{ \mathcal{X}_{K_1+1}, \dots, \mathcal{X}_{2K_1} \}$, $\dots$,  where each group has $K_1$ tensors in order. Then perform the tensor concatenations over these $K_2K_3\cdots K_N$ groups to obtain $K_2K_3\cdots K_N$ 1-mode row block tensors denoted by $\mathcal{X}^{(1)}_1 = \begin{bmatrix}
			\mathcal{X}_1 & \cdots & \mathcal{X}_{K_1}
		\end{bmatrix}_1, \dots, \mathcal{X}^{(1)}_{K_2K_3\cdots K_N} = \begin{bmatrix}
			\mathcal{X}_{|\mathbf{K}|-K_1+1} & \cdots & \mathcal{X}_{|\mathbf{K}|}
		\end{bmatrix}_1$, respectively, which are $2N$th-order paired tensors in $\mathbb{C}^{I_1\times J_1K_1 \times I_2 \times J_2 \times \cdots \times I_N \times J_N}$.
		
		\item Divide the  tensors $ \{ \mathcal{X}_1^{(1)}, \mathcal{X}_2^{(1)}, \dots,\mathcal{X}_{K_2\cdots K_N}^{(1)} \}$ obtained above sequentially into $K_3K_4 \cdots K_N$ groups $\{ \mathcal{X}_1^{(1)}, \dots, \mathcal{X}_{K_2}^{(1)} \}$, $\{ \mathcal{X}_{K_2+1}^{(1)},\dots, \mathcal{X}_{2K_2}^{(1)} \}$, $\dots$, where each group has $K_2$ tensors in order. Then perform the tensor concatenations over these $K_3K_4\cdots K_N$ groups of tensors to obtain $K_3K_4\cdots K_N$ $2$-mode row block tensors denoted by $\mathcal{X}_1^{(2)},  \dots, \mathcal{X}_{K_3K_4\cdots K_N}^{(2)} \in \mathbb{C}^{I_1 \times J_1 K_1 \times I_2 \times J_2 K_2 \times \cdots \times I_N \times J_N}$, respectively.
		
		\item  Keep repeating the above processes until the last $N$-mode row block tensor is obtained, which is the mode row block tensor $\mathcal{X}\in \mathbb{C}^{I_1\times J_1K_1\times \cdots \times I_N\times J_NK_N}$, denoted by $\mathcal{X} = \begin{vmatrix}
			\mathcal{X}_1 & \mathcal{X}_2 & \cdots  & \mathcal{X}_{|\mathbf{K}|}
		\end{vmatrix}$.
	\end{enumerate}
	Here the symbol $| \cdots |$ is used to distinguish it from the $n$-mode block tensors. The mode column block tensor can be defined in a similar way and is denoted by 
	\[\begin{vmatrix}
		\mathcal{X}_1^{\top} & \mathcal{X}_2^{\top}    & \cdots  & \mathcal{X}_{|\mathbf{K}|}^{\top}
	\end{vmatrix}^{\top} \in \mathbb{C}^{I_1K_1\times J_1 \times \cdots \times I_NK_N \times J_N}.\]
\end{definition}

\begin{definition}[Tensor blockings \cite{ragnarsson2012block}]
\label{Def: tensorblocking}
    For a tensor $\mathcal{A} \in \mathbb{C}^{L_1\times \cdots \times L_N}$, we say that $ \{ \mathbf{m}^{(1)}, \cdots , \mathbf{m}^{(N)} \}$ specifies  a blocking for $\mathcal{A}$, if for all $n = 1, \dots, N$, $\mathbf{m}^{(n)} = ( m_1^{(n)}, \dots, m_{I_n}^{(n)} )$ is an $I_n$-dimensional vector of positive integers that sum to $L_n$, which divides the $n$th mode of $\mathcal{A}$ into $I_n$ parts. The blocking identifies $\mathcal{A}$ as an $I_1 \times \cdots \times I_N$ block tensor, whose $\mathbf{i} = ( i_1, \cdots , i_N )$-th block is a subtensor of size $m_{i_1}^{(1)} \times \cdots \times m_{i_N}^{(N)}$, denoted by  $\mathcal{A}_{[\mathbf{i}]} := \mathcal{A}(l_{i_1-1}^{(1)}+1: l_{i_1}^{(1)}, \dots, l_{i_N-1}^{(N)}+1: l_{i_N}^{(N)})$, where $l^{(n)}_{i_n}$ denotes the sum of the first $i_n$ terms of the vector $\mathbf{m}^{(n)}$, that is, $l^{(n)}_{i_n} = {\textstyle \sum_{j = 1}^{i_n}} m_j^{(n)}$ for all $i_n=1,\dots,  I_n$ and $n=1,\dots,N$.
\end{definition}
\begin{remark}
	For example, the mode column block tensor $\left[\begin{smallmatrix}
		\mathcal{X}_1^{\top}  & \cdots  & \mathcal{X}_{|\mathbf{K}|}^{\top}
	\end{smallmatrix}\right]^{\top}$ in \cref{Def: modeblocktensor} has a blocking $\{ \mathbf{m}^{(1)},\dots, \mathbf{m}^{(2N)} \}$ with $\mathbf{m}^{(2n-1)} = ( \underbrace{ I_n,\dots, I_n }_{K_n \ \mathrm{terms}})$ and $\mathbf{m}^{(2n)} = J_n$ for $n=1,\dots, N$, and its subblocks are exactly the tensors $\mathcal{X}_m\ (m=1,\dots,|\mathbf{K}|)$ via the  blocking in \cref{Def: tensorblocking}.
\end{remark}

\subsection{Tensor U-eigenvalues and norms}

The tensor eigenvalue problem via the  Einstein product often arises in elastic mechanics \cite{mehrabadi1990eigentensors, itskov2000theory}, and is further studied in \cite{Brazell2013, cui2016eigenvalue, wang2022generalized}, among others. Chen et al. \cite{chen2019multilinear, chen2021multilinear} generalized the matrix-based Rayleigh quotient iteration method for computing tensor U-eigenvalues.  

\begin{definition} \label{def:U eig}(Tensor U-eigenvalue \cite{chen2021multilinear})
	Given $\mathcal{A}\in \mathbb{C}^{I_1\times I_1 \times \cdots \times I_N\times I_N}$, if an $N$th-order nonzero tensor $ \mathcal{X}\in \mathbb{C}^{I_1\times \cdots \times I_N}$ and $\lambda \in \mathbb{C}$ satisfy $\mathcal{A}*\mathcal{X}=\lambda \mathcal{X}$, 
	then $\lambda$ and $\mathcal{X}$ are called the U-eigenvalue and U-eigentensor of $\mathcal{A}$, respectively. The tensor $\mathcal{A}$ is said to be stable if all its U-eigenvalues are in the open left-half  plane.
\end{definition} 

\begin{remark} \label{rem:stability}
It can be easily shown that $\lambda$ is  a $U$-eigenvalue of a tensor $\mathcal{A}$ if and only if it is an eigenvalue of the associated matrix $\phi(\mathcal{A})$. Thus, the tensor $\mathcal{A}$  is  stable if and only if the  matrix $\phi(\mathcal{A})$ is Hurwitz stable.
\end{remark}

\begin{remark}
Another type of U-eigenvalues for tensors was defined and studied in \cite{NQB14, ZNZ20, CPZC21}, which should not be confused with the one used in this paper. 
\end{remark}

The tensor Schur decomposition for even-order paired tensors follows directly from \cite[Theorem 4.9]{liang2019tensor}. Thus, its proof is omitted.
\begin{lemma}[Tensor Schur decomposition]  
\label{lem:schur}
	Given $\mathcal{A} \in \mathbb{C}^{I_1 \times I_1 \times \cdots \times I_N \times I_N}$,   there exist a unitary tensor $\mathcal{U}$ and an upper triangular tensor $\mathcal{T}$ of the same size such that $\mathcal{A}=\mathcal{U} * \mathcal{T} * \mathcal{U}^{\mathrm{H}}$.	Moreover, the diagonal entries $\mathcal{T}_{i_1i_1 \ldots i_N i_N}$ of $\mathcal{T}$ are the U-eigenvalues of $\mathcal{A}$.
\end{lemma}

The tensor SVD for even-order tensors has been proposed in \cite{Brazell2013}, which for even-order {\it paired} tensors is given below.
\begin{lemma}[Tensor SVD]
 \label{Lem: tensorSVD}
		\it The tensor $\mathrm{SVD}$ of an even-order paired tensor $\mathcal{A}\in \mathbb{C}^{I_1\times J_1 \times \cdots \times I_N\times J_N }$ is
		 $\mathcal{A} = \mathcal{U} * \mathcal{D} * \mathcal{V}^{\rm H}$,  
		where $\mathcal{U} \in \mathbb{C}^{I_1\times I_1 \times \cdots \times I_N \times I_N}$ and $\mathcal{V} \in \mathbb{C}^{J_1\times J_1 \times \cdots \times J_N \times J_N}$ are unitary tensors, and $\mathcal{D} \in \mathbb{C}^{I_1\times J_1 \times \cdots \times I_N\times J_N }$ is a tensor with $\mathcal{D}_{i_1j_1\cdots i_Nj_N}  = 0$ for all $\mathbf{i} \ne \mathbf{j}$, whose diagonal entries $\mathcal{D}_{i_1 i_1\cdots i_Ni_N}$ are called the singular values of $\mathcal{A}$.   
\end{lemma} 

The outer product of two tensors $\mathcal{X}\in \mathbb{C}^{I_1\times \cdots \times I_N}$ and $\mathcal{Y}\in \mathbb{C}^{J_1\times \cdots \times J_M}$ is an $(N+M)$th-order tensor $\left(\mathcal{X}\circ \mathcal{Y}\right)_{i_1\cdots i_Nj_1\cdots j_M} = \mathcal{X}_{i_1\cdots i_N}\mathcal{Y}_{j_1\cdots j_M}$. The inner product of two  tensors $\mathcal{X}, \mathcal{Y} \in \mathbb{C}^{J_1\times \cdots \times J_N}$ is $\left \langle \mathcal{X}, \mathcal{Y} \right \rangle = \sum_{ \mathbf{1} \le \mathbf{j} \le \mathbf{J}} \overline{\mathcal{X}_{j_1\cdots j_N}}\mathcal{Y}_{j_1\cdots j_N}$, which gives the Frobenius norm
 $\left \| \mathcal{X} \right \|_F := \sqrt{\left \langle \mathcal{X}, \mathcal{X} \right \rangle } =  \sqrt{ \sum_{\mathbf{1} \le \mathbf{j} \le {\mathbf{J}}} \left | \mathcal{X}_{j_1\cdots j_N} \right | ^2 }$.  

\begin{definition}[Spectral norm \cite{Ma2019perturbation}] \label{def:spectrum} 
	The spectral norm of $\mathcal{A} $ is defined by $\left\| \mathcal{A} \right\|_2 : = \sqrt{\lambda_{\rm{max}}(\mathcal{A}^{\mathrm{H}} * \mathcal{A})}$,	where $\lambda_{\rm{max}}(\cdot)$ denotes the largest U-eigenvalue.
\end{definition}

\begin{remark}[\cite{wang2022generalized}] \label{rem:sigma}
	The spectral norm of a tensor $\mathcal{A}$ satisfies $\| \mathcal{A} \|_2 = \sigma_{\rm max}(\mathcal{A})$, namely the largest singular value of $\mathcal{A}$.  Moreover, it is compatible with the Frobenius norm, that is, $\left \| \mathcal{A}* \mathcal{X} \right \|_F \le \left\| \mathcal{A} \right\|_2 \left \| \mathcal{X} \right \|_F$ for all $\mathcal{X} \in \mathbb{C}^{J_1\times \cdots \times J_N}$.
\end{remark}

\section{Hamiltonian and symplectic tensors, and tensor Kronecker product}
\label{sec. 3}

\subsection{Hamiltonian and  symplectic tensors}

Ragnarsson and Van Loan \cite{ragnarsson2012block} studied the unfolding patterns of block tensors: the subblocks of a tensor can be mapped to contiguous blocks in the unfolding matrix through a series of row and column permutations. Specifically, let $s = qr$.  A perfect shuffle permutation $\Pi _{q,r} \in \mathbb{R}^{s\times s}$ is defined by
\[ \Pi_{q,r} \mathbf{z} = \begin{bmatrix}
	z_{1:r:s} \\
	z_{2:r:s} \\
	\vdots  \\
	z_{r:r:s}
\end{bmatrix}, \ \ \forall  \ \mathbf{z}\in \mathbb{C}^{s}. \]

 The following lemma is an immediate consequence of \cite[Theorem 3.3]{ragnarsson2012block}. 

\begin{lemma} \label{Lemma: unfolding}
	Given even-order paired tensors $\mathcal{A}, \mathcal{B}, \mathcal{C}, \mathcal{D} \in \mathbb{C}^{I_1\times I_1 \times \cdots \times I_N \times I_N}$, there exists a permutation matrix $P = Q_N\cdots Q_2Q_1$ such that 
	\[ \phi\left( \begin{bmatrix}
		\mathcal{A} & \mathcal{B} \\
		\mathcal{C} & \mathcal{D}
	\end{bmatrix}_n \right) =  P \begin{bmatrix}
		\phi(\mathcal{A})  & \phi(\mathcal{B}) \\
		\phi(\mathcal{C}) & \phi(\mathcal{D})
	\end{bmatrix} P^{\top}, \]
	where $Q_k = I_{2I_1\cdots I_N}$ for $k \le n$, and $Q_k = I_{I_{k+1}\cdots I_N}\otimes \Pi _{I_k,2} \otimes I_{I_1\cdots I_{k-1}}$ for $k \ge n+1$.
\end{lemma}

The results to be derived in the sequel can be done  for all $n$-mode block tensors for $n=1,\dots,N$. For ease of representation, we present them for the $1$-mode block tensor case and also omit the subscript, i.e., let $ \left[\begin{smallmatrix}
	\mathcal{A} & \mathcal{B} \\
	\mathcal{C} & \mathcal{D}
\end{smallmatrix}\right] = \left[\begin{smallmatrix}
	\mathcal{A} & \mathcal{B} \\
	\mathcal{C} & \mathcal{D}
\end{smallmatrix}\right]_1 \in \mathbb{C}^{2I_1\times 2I_1\times I_2\times I_2 \times \cdots \times I_N\times I_N}$, for $\mathcal{A}, \mathcal{B}, \mathcal{C}, \mathcal{D}\in \mathbb{C}^{I_1\times I_1 \times \cdots \times I_N\times I_N}$ and $\left[\begin{smallmatrix}
	\mathcal{X} \\ \mathcal{Y}
\end{smallmatrix}\right] =\left[\begin{smallmatrix}
	\mathcal{X} \\
	\mathcal{Y}
 \end{smallmatrix}\right]_1$ for $N$th-order tensors $\mathcal{X}, \mathcal{Y}\in \mathbb{C}^{I_1\times \cdots \times I_N}$. We introduce the block-structured tensor $		\mathcal{J} = \left[\begin{smallmatrix}
	\mathcal{O} & \mathcal{I} \\
	-\mathcal{I} & \mathcal{O}
		\end{smallmatrix}\right]$,
where $\mathcal{I}, \mathcal{O} $ are the identity tensor and zero tensor, respectively. It reduces to the well-known symplectic matrix $J=\left[\begin{smallmatrix}
	O & I \\
	-I & O
\end{smallmatrix}\right]$ if $N = 1$.

\begin{definition}[Symplectic tensor]\label{def:sym}
	We call $\mathcal{S}\in \mathbb{C}^{2I_1\times 2I_1 \times I_2 \times I_2 \times \cdots \times I_N \times I_N}$ a symplectic tensor if $\mathcal{S}^{\mathrm{H}}* \mathcal{J}* \mathcal{S} = \mathcal{J}$.
\end{definition}

\begin{definition} [Hamiltonian tensor]\label{def:ham}
	We call $\mathcal{M} \in \mathbb{C}^{2I_1\times 2I_1 \times I_2\times I_2 \times \cdots \times I_N\times I_N}$ a Hamiltonian tensor if $(\mathcal{J}*\mathcal{M})^{\mathrm{H}} = \mathcal{J}* \mathcal{M}$.
\end{definition}
These notions reduce to symplectic and Hamiltonian matrices when $N = 1$.

\begin{remark}
It is worthwhile to point out that Hamiltonian and symplectic structures of tensors are not preserved under the isomorphic map $\phi$.  
\end{remark}

The following proposition follows directly from \cref{Pro: blocktensor} and the definition of Hamiltonian tensors. Thus, its proof is omitted.
\begin{proposition}
	 $\mathcal{M} \in \mathbb{C}^{2I_1\times 2I_1 \times I_2\times I_2 \times \cdots \times I_N\times I_N}$ is a Hamiltonian tensor if and only if it  has the block structure $\mathcal{M} = \left[\begin{smallmatrix}
			\mathcal{A} & \mathcal{G} \\
			\mathcal{K} & -\mathcal{A}^{\mathrm{H}}
		\end{smallmatrix}\right]$, where $\mathcal{A}, \mathcal{G}, \mathcal{K} \in \mathbb{C}^{I_1\times I_1 \times \cdots \times I_N  \times I_N}$ with  $\mathcal{G}, \mathcal{K}$ being Hermitian.
\end{proposition}

Analogous to the matrix case,  Hamiltonian tensors enjoy the following spectral properties under the tensor algebra via the  Einstein product.  The proof is omitted.
\begin{proposition}
\label{prop: Hamiltonian}
	 Let $\lambda$ be a U-eigenvalue and the nonzero $N$th-order tensor $\left[\begin{smallmatrix}
		\mathcal{X} \\
		\mathcal{Y}
	\end{smallmatrix}\right] \in \mathbb{C}^{2I_1\times I_2 \times \cdots \times I_N}$ be its associated U-eigentensor of a Hamiltonian tensor $\mathcal{M}$. Then
	\begin{enumerate}
		\rm \item \it $ -\bar{\lambda}$ is also a U-eigenvalue of $\mathcal{M}$; and $\bar{\lambda}, -\lambda$ are U-eigenvalues too if $\mathcal{M}$ is real. 

		\rm \item \it $\mathcal{X}^{\mathrm{H}}* \mathcal{K}* \mathcal{X} + \mathcal{Y}^{\mathrm{H}}* \mathcal{G}* \mathcal{Y} = \left( \lambda+\bar{\lambda} \right) \mathcal{X}^{\mathrm{H}}* \mathcal{Y}$;
				
		\rm \item \it  $\mathcal{G}, \mathcal{K}$ are positive definite $\Rightarrow \mathrm{Re}(\lambda) \ne 0$;
		
		\rm \item \it $\mathcal{G}, \mathcal{K}$ are positive semidefinite and $\mathrm{Re}(\lambda)=0\Rightarrow \mathcal{A}*\mathcal{X} = \lambda \mathcal{X}, \mathcal{A}^{\mathrm{H}}*\mathcal{Y} = -\lambda\mathcal{Y}$.
	\end{enumerate}
\end{proposition}

\begin{theorem} [Schur-Hamiltonian tensor decomposition] 
\label{thm:Schur_Hamiltonian}
	Let $\mathcal{M}$ be a Hamiltonian tensor, and have no U-eigenvalues on the imaginary axis. Then there exists a unitary symplectic tensor $\mathcal{Q}$ such that
\begin{equation} \label{eq:Schur-Hamiltonian}
     \mathcal{Q}^{\mathrm{H}}* \mathcal{M}* \mathcal{Q} = \begin{bmatrix}
		\mathcal{T}  & \mathcal{R} \\
		\mathcal{O}  & -\mathcal{T}^{\mathrm{H}}
	\end{bmatrix},  
\end{equation}	
	where $\mathcal{T}$ is an upper triangular tensor whose U-eigenvalues are in the open left-half plane and $\mathcal{R}$ is a Hermitian tensor.
\end{theorem}
\begin{proof}
	We carry out the proof by using tensor unfolding. Let matrices $A = \phi(\mathcal{A}), G = \phi(\mathcal{G}), K = \phi(\mathcal{K})$. By \cref{Lemma: unfolding}, there is a permutation matrix $P$ such that
    \[ P^{\top}\phi(\mathcal{M})P =  \begin{bmatrix}
			A & G \\
			K & -A^{\mathrm{H}}
		\end{bmatrix} =: M,\ \mathrm{and}\
  P^{\top}\phi(\mathcal{J})P =  \begin{bmatrix}
			O & I \\
			-I & O
		\end{bmatrix}. \]
As $\mathcal{M}$  has no U-eigenvalues on the imaginary axis, the matrix $M$  obtained above has no eigenvalues on the imaginary axis. According to \cite[Theorem 3.1]{paige1981schur}, we perform the Schur-Hamiltonian decomposition for the Hamiltonian matrix $M$, that is, there exists a unitary symplectic matrix $Q$ such that $Q^{\mathrm{H}} M Q = \left[\begin{smallmatrix}
		T & R \\
		O & -T^{\mathrm{H}}
	\end{smallmatrix}\right]$, where $T,R \in \mathbb{C}^{\left | \mathbf{I} \right | \times \left | \mathbf{I} \right | },$	 $T$ is upper triangular whose eigenvalues are in the open left half complex plane, and $R^{\mathrm{H}} = R$. Let $\mathcal{T} = \phi^{-1}(T), \mathcal{R} = \phi^{-1}(R)$. Then we have
	\begin{align} \label{eq:30jan}
PQ^{\mathrm{H}}P^{\top}\phi(\mathcal{M})PQP^{\top}  = P \begin{bmatrix}
		T & R \\
		O & -T^{\mathrm{H}}
	\end{bmatrix} P^{\top} 
 	= \phi\left( \begin{bmatrix}
		\mathcal{T}  & \mathcal{R} \\
		\mathcal{O}  & -\mathcal{T}^{\mathrm{H}}
	\end{bmatrix} \right).
	\end{align}
	Define $\mathcal{Q} := \phi^{-1}(PQP^{\top})$. It follows from \cref{eq:30jan} and the isomorphic property of $\phi^{-1}$ we get \cref{eq:Schur-Hamiltonian}. Clearly, $\mathcal{R}$ is a Hermitian tensor. Next, we show that $\mathcal{Q}$ is a unitary symplectic tensor.  Notice that
	\begin{align*}
		\mathcal{Q}^{\rm H}* \mathcal{Q} & = \phi^{-1} (PQ^{\mathrm{H}}P^{\top})* \phi^{-1}(PQP^{\top})  = \phi^{-1}(PQ^{\mathrm{H}}P^{\top}PQP^{\top}) = \phi^{-1}(I) = \mathcal{I}.
	\end{align*}
 	Thus, $\mathcal{Q}$ is unitary. Moreover, observe that
 	\begin{align*}
\phi(\mathcal{Q})^{\mathrm{H}}\phi(\mathcal{J})\phi(\mathcal{Q}) &= (PQP^{\top})^{\mathrm{H}}(PJP^{\top})(PQP^{\top}) \\
 		& = PQ^{\mathrm{H}}JQP^{\top} = PJP^{\top} = \phi(\mathcal{J}).
 	\end{align*}
	Applying $\phi^{-1}$ on both sides of the above equation gives that $\mathcal{Q}^{\mathrm{H}}* \mathcal{J} * \mathcal{Q} = \mathcal{J}$. Hence $\mathcal{Q}$ is a symplectic tensor.
\end{proof} 

The following proposition can be derived from the definitions of unitary tensors and symplectic tensors and \cref{Pro: blocktensor}, thus its proof is omitted.
\begin{proposition}
	\label{Pro: Symblock}
	$\mathcal{S} \in \mathbb{C}^{2I_1\times 2I_1 \times I_2 \times I_2 \times \cdots \times I_N\times I_N}$ is a unitary symplectic tensor if and only if $\mathcal{S}$ has the block structure $\mathcal{S} = \left[\begin{smallmatrix}
		\mathcal{Q}_1 & \mathcal{Q}_2 \\
		-\mathcal{Q}_2 & \mathcal{Q}_1
	\end{smallmatrix}\right]$, where $\mathcal{Q}_1, \mathcal{Q}_2 \in \mathbb{C}^{I_1\times I_1 \times \cdots \times I_N  \times I_N}$ satisfy that $\mathcal{Q}_1^{\rm H}*\mathcal{Q}_2$ is  Hermitian and $\mathcal{Q}_1^{\rm H}*\mathcal{Q}_1 + \mathcal{Q}_2^{\rm H}*\mathcal{Q}_2 = \mathcal{I}$.
\end{proposition}

Paige and Van Loan \cite{paige1981schur} proved that a special SVD exists for unitary symplectic matrices. Based on \cref{Pro: Symblock} we can generalize this result to unitary symplectic tensors, whose proof is omitted due to page limit.
\begin{theorem}[Symplectic tensor SVD]
	\label{Thm: SymSVD}
	If $\mathcal{Q} \in \mathbb{C}^{2I_1  \times 2I_1 \times I_2 \times I_2 \times \cdots \times I_N \times I_N }$ is a unitary symplectic tensor, then there exist unitary tensors $\mathcal{U}, \mathcal{V} \in \mathbb{C}^{I_1\times I_1 \times \cdots \times I_N \times I_N }$ such that $\left[\begin{smallmatrix}
		\mathcal{U}^{\mathrm{H}}  & \mathcal{O} \\
		\mathcal{O} & \mathcal{U}^{\mathrm{H}}
	\end{smallmatrix}\right]* \mathcal{Q} * \left[\begin{smallmatrix}
		\mathcal{V}  & \mathcal{O} \\
		\mathcal{O} & \mathcal{V}
	\end{smallmatrix}\right] = \left[\begin{smallmatrix}
		\mathcal{S} & \mathcal{D} \\
		-\mathcal{D}	& \mathcal{S}
	\end{smallmatrix}\right]$, where $\mathcal{S},\mathcal{D} \in \mathbb{C}^{I_1\times I_1 \times \cdots \times I_N \times I_N}$ are diagonal tensors with $0\le \mathcal{S}_{i_1i_1\cdots i_Ni_N} \le 1$
	and $(\mathcal{S}_{i_1i_1\cdots i_Ni_N})^2 + (\mathcal{D}_{i_1i_1\cdots i_Ni_N})^2 = 1$.
\end{theorem}

\subsection{Tensor Kronecker product and Vec operation}
\label{Section: 3.2}
The Kronecker product and vectorized operation are frequently used in matrix computations \cite{graham2018kronecker}. In this section, we generalize these operations to even-order paired tensors.

A definition of tensor Kronecker product is introduced \cite{pickard2023kronecker}. Here we give this definition by using the tensor blocking technique (\cref{Def: tensorblocking}).
\begin{definition}[Tensor Kronecker product]
	\label{TensorKronecker}
	Given tensors $\mathcal{A}\in \mathbb{C}^{K_1\times \cdots \times K_N}$ and $\mathcal{B} \in \mathbb{C}^{L_1\times \cdots \times L_N}$, we define their Kronecker product $\mathcal{A} \otimes \mathcal{B} \in \mathbb{C}^{K_1L_1\times \cdots \times K_NL_N}$ as a $K_1\times \cdots \times K_N$ block tensor having the blocking $\left \{ \mathbf{m}^{(1)}, \cdots, \mathbf{m}^{(N)} \right \}$, where
	\begin{align*}
		\mathbf{m}^{(n)} = 
		(\underbrace{L_n,\cdots, L_n}  _{K_n\ \mathrm{terms}}) ,\ 1\le n \le N. 
	\end{align*}
	There are $|\mathbf{K} | = K_1\cdots K_N$ subblocks in total, and the $\mathbf{k}$-th subblock is the $L_1\times \cdots \times  L_N$ tensor $\mathcal{A}_{k_1\cdots k_N} \mathcal{B}$, that is,
	$(\mathcal{A}\otimes \mathcal{B})_{[\mathbf{k}]}
		= \mathcal{A}_{k_1\cdots k_N} \mathcal{B}$ for all $\mathbf{1}\le \mathbf{k}\le \mathbf{K}$.
\end{definition}
\begin{remark}
	We can also perform the tensor Kronecker product between two tensors with different orders, after supplementing the mode dimensions by 1 for the  tensor of smaller order.
\end{remark}

The following result tells us that the tensor Kronecker product is an extension of the matrix Kronecker product. 

\begin{proposition}
	\label{Pro: TensorKron}
	Let $\mathcal{A}\in \mathbb{C}^{I_1\times J_1 \times \cdots \times I_N \times J_N}, \mathcal{B} \in \mathbb{C}^{K_1\times L_1 \times \cdots \times K_N \times L_N}$ be even-order paired tensors. Then
	\begin{enumerate}
        \item $rank_U(\mathcal{A}\otimes \mathcal{B}) = rank_U(\mathcal{A})\cdot rank_U(\mathcal{B})$.
		\item $(\mathcal{A}\otimes \mathcal{B})^{\mathrm{H}} = \mathcal{A}^{\mathrm{H}} \otimes \mathcal{B}^{\mathrm{H}}$.
        \item $(\mathcal{A}\otimes \mathcal{B}) * (\mathcal{C}\otimes \mathcal{D}) = (\mathcal{A}*\mathcal{C})\otimes (\mathcal{B}* \mathcal{D})$, for $\mathcal{C}\in \mathbb{C}^{J_1\times R_1 \times \cdots \times J_N \times R_N}$ and  $\mathcal{D}\in \mathbb{C}^{L_1\times S_1 \times \cdots \times L_N \times S_N}$.
    \item If $\mathcal{A}, \mathcal{B}$ are invertible, then $\mathcal{A}\otimes\mathcal{B}$ is  invertible, and $\left( \mathcal{A}\otimes \mathcal{B} \right)^{-1} = \mathcal{A}^{-1} \otimes \mathcal{B}^{-1}$.
	\end{enumerate}
\end{proposition}

\begin{proof}
	The proof is given in Appendix B.
\end{proof}

\begin{corollary}
	\label{Pro: tensorKronQR}
    For the Kronecker product tensor $\mathcal{A} = \mathcal{B}\otimes \mathcal{C}$, its tensor decomposition via the Einstein product can be obtained from the tensor decompositions of $\mathcal{B}$ and $\mathcal{C}$.  Specifically, $\mathcal{A} = \mathcal{U}* \mathcal{T}*\mathcal{U}^{\rm H}$ with $\mathcal{U}=\mathcal{U}_1\otimes \mathcal{U}_2$ and $\mathcal{T}=\mathcal{T}_1\otimes \mathcal{T}_2$ is the tensor Schur decomposition of $\mathcal{A}$, if we have the tensor Schur decompositions  $\mathcal{B} = \mathcal{U}_1 * \mathcal{T}_1*\mathcal{U}_1^{\rm H}$ and $\mathcal{C} = \mathcal{U}_2 * \mathcal{T}_2 * \mathcal{U}_2^{\rm H}$. As a result, the U-eigenvalues of $\mathcal{A}$ are the products of all U-eigenvalues of $\mathcal{B}$ and $\mathcal{C}$.
\end{corollary}
\begin{proof}
	The result follows from \cref{Pro: TensorKron} and the fact that the Kronecker product of two upper-triangular tensors is still upper-triangular.
\end{proof}

The vectorized operation `Vec' for a matrix $A \in \mathbb{C}^{m\times n}$ stacks its columns sequentially, i.e., Vec$(A) = \left[ A_{:, 1}^{\top}, \cdots , A_{:, n}^{\top} \right]^{\top} \in \mathbb{C}^{nm}$ \cite{golub2013matrix}. Next, we define the Vec operation for even-order paired tensors that transforms a $2N$th-order paired tensor to an $N$th-order tensor based on \cref{Def: modeblocktensor} of mode block tensors, which reduces to the matrix vectorization  when $N = 1$.
\begin{definition}[Tensor Vec operation]
	\label{Def: tensorVec}
	Given an even-order paired tensor $\mathcal{X}\in \mathbb{C} ^{J_1 \times K_1 \times \cdots \times J_N  \times K_N}$, let $\mathcal{X}_{\mathrm{ivec}(\mathbf{k}, \mathbf{K})} = \mathcal{X}_{:k_1, \cdots, :k_N} \in \mathbb{C}^{J_1 \times \cdots \times J_N}$ be its  subtensor for each given set of the mode column indices $\mathbf{k}=(k_1, \ldots, k_N)$. We define $\mathrm{Vec}(\mathcal{X}) \in \mathbb{C}^{K_1J_1 \times \cdots \times K_NJ_N}$ as the mode column block tensor $\mathrm{Vec}(\mathcal{X}):=\begin{vmatrix}
		\mathcal{X}_1^{\top}  & \mathcal{X}_2^{\top} & \cdots & \mathcal{X}_{|\mathbf{K}|}^{\top}
	\end{vmatrix}^{\top}$.
\end{definition}
    
    Equivalently, the tensor Vec operation can be defined using  tensor blockings (\cref{Def: tensorblocking}), as shown below.
\begin{proposition} \label{prop: TVO1}
	Let $\mathcal{X}\in \mathbb{C} ^{J_1 \times K_1 \times \cdots \times J_N \times K_N}$ be an even-order paired tensor. Then $\mathrm{Vec}(\mathcal{X}) \in \mathbb{C}^{K_1J_1 \times \cdots \times K_NJ_N}$ is an $N$th-order block tensor with $| \mathbf{K} |$ subblocks, where each subblock is a tensor of size $J_1 \times \cdots \times J_N$. Concretely, it has the blocking $\{ \mathbf{m}^{(1)}, \cdots, \mathbf{m}^{(N)} \}$, where 
    \[ \mathbf{m}^{(n)} = (\underbrace{ J_n,\cdots, J_n  }_{K_n\ \mathrm{terms}}) ,\ 1\le n \le N, \]
    and the $ \mathbf{k} $th subblock tensor is $\mathrm{Vec}(\mathcal{X})_{[\mathbf{k}]} = \mathcal{X}_{:k_1, \cdots,:k_N}$.
\end{proposition}
\begin{proof}
	The result directly follows  from the \cref{Def: tensorblocking} of tensor blocking and \cref{Def: tensorVec} of tensor $\mathrm{Vec}$ operation.
\end{proof}

If tensors are given in the generalized CPD format (see \cref{Def: GCPD} in Appendix A), their tensor Kronecker product and Vec operation can be computed in an economic way by keeping the original structured format. In other words, the resultant tensor of these operations can be expressed from the Kronecker product and vectorization of the factor matrices. More results on fast computational methods for the structured tensors are given in \cref{appendix C}.

\begin{proposition} \label{prop: TVO2}
	Given two $2N$th-order paired tensors in the generalized CPD format $\mathcal{A} = \sum_{r=1}^{R} A_1^{(r)}\circ A_2^{(r)}\circ \cdots \circ A_N^{(r)} $ and $\mathcal{B} = \sum_{s=1}^{S} B_1^{(s)}\circ B_2^{(s)}\circ \cdots \circ B_N^{(s)}$, we have
	\begin{enumerate}
		\item $\mathcal{A} \otimes \mathcal{B} = \sum_{r=1}^R \sum_{s=1}^S (A_1^{(r)} \otimes B_1^{(s)}) \circ \cdots \circ (A_N^{(r)} \otimes B_N^{(s)}) $;
		\item $\mathrm{Vec}(\mathcal{A}) = \sum_{r=1}^R \mathrm{Vec}(A_1^{(r)}) \circ \cdots \circ \mathrm{Vec}(A_N^{(r)})$. 
	\end{enumerate}
\end{proposition}
\begin{proof}
	The result follows from the definitions of the tensor Kronecker product, Vec operation, and outer product.
\end{proof}

Next, we derive a very useful result, which reduces to  the matrix case  when $N=1$.
\begin{proposition}	
	\label{Pro: KronVec}
 Given $\mathcal{U} \in \mathbb{C}^{I_1\times J_1 \times \cdots \times I_N \times J_N}$, $\mathcal{X}\in \mathbb{C}^{J_1\times K_1 \times \cdots \times J_N \times K_N}$, $\mathcal{W}\in \mathbb{C}^{K_1\times L_1 \times \cdots}$ $^{\times K_N \times L_N}$, we have $	\mathrm{Vec}(\mathcal{U}* \mathcal{X}* \mathcal{W}) = (\mathcal{W}^{\top} \otimes \mathcal{U})* \mathrm{Vec}(\mathcal{X})$.
\end{proposition}
\begin{proof}
	The proof is given in \cref{appemd:proofs}. 
\end{proof}

\begin{lemma}
	\label{Lem: KronUeigen}
	Let $\{ \lambda_1, \cdots, \lambda_{|\mathbf{I}|} \}$ and $\{ \mu_1, \cdots , \mu_{|\mathbf{J}|} \}$ be the U-eigenvalues of $\mathcal{A} \in \mathbb{C}^{I_1\times I_1 \times \cdots \times I_N \times I_N}$ and $\mathcal{B} \in \mathbb{C}^{J_1 \times J_1 \times \cdots \times J_N \times J_N}$, respectively. Then the U-eigenvalues of $\left(\mathcal{I}_{\mathbf{J}} \otimes \mathcal{A} - \mathcal{B}^{\top}\otimes \mathcal{I}_{\mathbf{I}} \right)$ are $\left\{ \lambda_i - \mu_j \big| \ i = 1,\cdots, |\mathbf{I}|; j = 1, \cdots, |\mathbf{J}| \right\}$.
\end{lemma}
\begin{proof}
	According to \cref{lem:schur}, we perform tensor Schur decomposition for $\mathcal{A}$ and $\mathcal{B}$ respectively, 
	 $\mathcal{A} = \mathcal{U}^{\mathrm{H}} * \mathcal{T}* \mathcal{U},\ \mathcal{B} ^{\top} = \mathcal{V}^{\mathrm{H}}* \mathcal{S}* \mathcal{V}$, where $\mathcal{U}$ and $\mathcal{V}$ are unitary, $\mathcal{T}$ and $\mathcal{S}$ are upper-triangular even-order paired tensors, whose diagonal elements are the U-eigenvalues of $\mathcal{A}$ and $\mathcal{B}$, respectively. From \cref{Pro: TensorKron}, we have
    \begin{align*}
    \mathcal{I}_{\mathbf{J}} \otimes \mathcal{A} - \mathcal{B}^{\top}\otimes \mathcal{I}_{\mathbf{I}} 
	= & \left( \mathcal{V}^{\mathrm{H}}* \mathcal{I}_{\mathbf{J}}* \mathcal{V} \right) \otimes \left( \mathcal{U}^{\mathrm{H}}* \mathcal{T}* \mathcal{U} \right) - \left( \mathcal{V}^{\mathrm{H}}* \mathcal{S} * \mathcal{V} \right)\otimes \left( \mathcal{U}^{\mathrm{H}}* \mathcal{I}_{\mathbf{I}}* \mathcal{U} \right) \\
	=& \left( \mathcal{V}\otimes \mathcal{U} \right)^{\mathrm{H}} * \left( \mathcal{I}_{\mathbf{J}}\otimes \mathcal{T} - \mathcal{S}\otimes \mathcal{I}_{\mathbf{I}} \right)* \left( \mathcal{V}\otimes \mathcal{U} \right).
	\end{align*}
	Noticing that $\left( \mathcal{V}\otimes \mathcal{U} \right)^{\mathrm{H}}* \left( \mathcal{V}\otimes \mathcal{U} \right) = \mathcal{I}_{\mathbf{J}}\otimes \mathcal{I}_{\mathbf{I}} = \mathcal{I}_{\mathbf{J}\odot \mathbf{I}}$, therefore $\mathcal{V}\otimes \mathcal{U}$ is a unitary tensor. Since $\mathcal{I}_{\mathbf{J}}\otimes \mathcal{T} - \mathcal{S}\otimes \mathcal{I}_{\mathbf{I}}$ has an upper-triangular structure, whose diagonal entries are $\{ \lambda_i - \mu_j \big| \ i = 1,\cdots, |\mathbf{I}|; j = 1, \cdots, |\mathbf{J}| \}$. From \cref{lem:schur}, the result follows. 
\end{proof}

Based on \cref{Pro: KronVec} and \cref{Lem: KronUeigen}, we obtain a necessary and sufficient condition for the existence and uniqueness of the solution of the Sylvester tensor equation  
\begin{equation} \label{Sylvester}
\mathcal{A}*\mathcal{E}+\mathcal{E}*\mathcal{B} = \mathcal{K},
\end{equation}
for which  numerical algorithms have been proposed in e.g. \cite{ali2016krylov, hajarian2020conjugate, huang2021numerical}.

\begin{theorem}
	\label{Thm: Sylvester}
	The Sylvester tensor equation \cref{Sylvester} has a unique solution if and only if $\mathcal{A}$ and $-\mathcal{B}$ have no common U-eigenvalues.
\end{theorem}	
\begin{proof}
    Performing the Vec operation on both sides of the Sylvester tensor equation \cref{Sylvester}, we get the following equivalent multilinear system $ \left( \mathcal{I} \otimes \mathcal{A} + \mathcal{B}^{\top}\otimes \mathcal{I} \right) * \mathrm{Vec}(\mathcal{E}) = \mathrm{Vec}(\mathcal{K})$. By \cref{Lem: KronUeigen}, the tensor $\mathcal{I} \otimes \mathcal{A} + \mathcal{B}^{\top}\otimes \mathcal{I}$ is invertible if and only if $\mathcal{A}$ and $-\mathcal{B}$ have no common U-eigenvalues. Therefore, the conclusion holds.
\end{proof}

Finally, we present the Lyapunov tensor equation under the Einstein product, which will be used to derive a tensor version of the bounded real lemma in \cref{thm:BRL}.
\begin{corollary}\label{lem: LTE}
	Let $\mathcal{A}\in \mathbb{C}^{I_1\times I_1\times \cdots \times I_N\times I_N}$ be stable and $\mathcal{Q}\in \mathbb{C}^{I_1\times I_1\times \cdots \times I_N\times I_N}$ be positive semidefinite. Then the Lyapunov tensor equation
	\begin{align} \label{Eq: LTE}
		\mathcal{A}^{\rm H}*\mathcal{E} + \mathcal{E}*\mathcal{A} + \mathcal{Q} = \mathcal{O} 
	\end{align} 
	has a unique solution, which is positive semidefinite.
\end{corollary}
\begin{proof}
	Since $\mathcal{A}$ is stable, $\mathcal{A}^{\rm H}$ and $-\mathcal{A}$ have no common U-eigenvalues. By \cref{Thm: Sylvester}, equation \cref{Eq: LTE} has a unique solution. Next, define the tensor exponential (\cite{elhalouy2022tensor}) $\mathrm{exp}(\mathcal{A}) := \sum_{k=0}^{\infty} \frac{\mathcal{A}^k}{k !}$,	where $\mathcal{A}^0=\mathcal{I}$ and $\mathcal{A}^k=\mathcal{A}* \mathcal{A}^{k-1},\ k\ge 1$.  Define $\mathcal{Y}(t) = \exp(t\mathcal{A}^{\rm H})*\mathcal{Q}*\exp(t\mathcal{A})$. Clearly, $\mathcal{Y}(0) = \mathcal{Q}$. By means of the isomorphic map $\phi$ defined in \cref{def:phi}, we can get a matrix $A=\phi(\mathcal{A})$. As  $\phi$ is an isomorphism, the matrix $A$ is stable as the tensor $\mathcal{A}$ is stable. Thus $\lim_{t\to+\infty} \mathrm{exp}(tA)=O$, which means $\lim_{t\to+\infty} \mathrm{exp}(t\mathcal{A})=\mathcal{O}$ as $\mathrm{exp}(t\mathcal{A}) = \phi^{-1}(\mathrm{exp}(tA))$. Consequently, $\lim_{t \to +\infty} \mathcal{Y}(t) = \mathcal{O}$. Noticing that $\dot{\mathcal{Y}}(t) = \mathcal{A}^{\rm H}*\mathcal{Y}(t) + \mathcal{Y}(t)*\mathcal{A}$, we have
\[
		\mathcal{Y}(t) - \mathcal{Y}(0) = \mathcal{A}^{\rm H}* \int_{0}^{t} \mathcal{Y}(s)ds + \int_{0}^{t} \mathcal{Y}(s)ds  *\mathcal{A}.
\]
	Sending  $t \to +\infty$ yields that  the positive semidefinite tensor $\int_{0}^{+\infty} \mathcal{Y}(s)ds$ is actually the solution of \cref{Eq: LTE}.
\end{proof}

\section{Algebraic Riccati tensor equation}
\label{sec. 4}
In this section, we introduce the basic control theory for the continuous-time MLTI systems \cref{MLTI-continuous}, and establish the unique existence of the positive semidefinite solution to the ARTE \cref{ARTE}, and prove the tensor version of the bounded real lemma and the small gain theorem. Moreover, a perturbation analysis is conducted for the ARTE \cref{ARTE}.

\subsection{Continuous-time multilinear time-invariant control systems}

The discrete-time MLTI system \cref{MLTI-Einstein} has been studied in \cite{chen2019multilinear,chen2021multilinear}. Fundamental control-theoretic notions, such as transfer function, stability, reachability and observability, are extended from discrete-time LTI systems to  discrete-time MLTI systems; see  Definition 3.4 for transfer function, Proposition  5.1 for stability, Definition 5.5 for reachability,  and Definition 5.13 for observability in Ref. \cite{chen2021multilinear} for details. These notions appear similar to their LTI counterparts. However, their criteria are established by means of tensor operations instead of matrix operations. 

In this section, we study the continuous-time MLTI system \cref{MLTI-continuous}. Similar to what has been done in Ref. \cite{chen2021multilinear} for the discrete-time case, the notions of transfer function, stability, controllability,  observability, stabilizability, and detectability can be extended to the continuous-time MLTI system \cref{MLTI-continuous} from the continuous-time LTI system theory. Like the discrete-time case,  these definitions appear similar to their LTI counterparts. For example, the continuous-time MLTI system \cref{MLTI-continuous} is stable if  all U-eigenvalues of the tensor $\mathcal{A}$ are in the open left-half plane.  Readers are referred to \cite{chen2021multilinear} to see the discrete-time counterparts of most of these definitions.

The transfer  tensor of the continuous-time MLTI system \cref{MLTI-continuous} is defined as  
\begin{equation} \label{eq:tf}
	\mathcal{G}(s) =    \mathcal{D} +  \mathcal{C}*(s \mathcal{I}- \mathcal{A})^{-1}*  \mathcal{B}. 
\end{equation}

Analogous to the matrix case, we say that $(\mathcal{A},\mathcal{B},\mathcal{C},\mathcal{D})$ is a state-space realization of the transfer tensor $\mathcal{G}(s)$. Similarly, we say that the state-space realization $(\mathcal{A},\mathcal{B},\mathcal{C},\mathcal{D})$ is  stable if  the tensor ${\cal A}$ is stable (see   \cref{def:U eig} and \cref{rem:stability}). Finally, the transfer tensor $\mathcal{G}(s)$ is said to be stable if all the roots of the denominator of each entry are in the open left-half plane; or equivalently, the unfolding determinant det$_U(\mathcal{G}(s))$ has no poles in the closed right-half plane.

 Recall that $\mathcal{RH}_\infty$ is the space of all proper and real-rational stable transfer matrices; see for example \cite[pp. 100]{ZDG96}.  As tensors  studied in this paper can be complex-valued, with a slight abuse of notation, we use $\mathcal{FH}_\infty$ to denote the space of all proper and rational stable transfer tensors. Here ``$\mathcal{F}$'' can be $\mathbb{R}$ or $\mathbb{C}$.  ``$\mathcal{F}$'' is used instead of ``$\mathcal{R}$'' as the coefficients of the  transfer tensors can be either real- or complex-valued.

Similar to the matrix case (\cite[pp. 100]{ZDG96}), for $\mathcal{G}(s) \in \mathcal{FH}_\infty$, we define its $H_{\infty}$  norm   as
\begin{equation} \label{eq:H infty norm}
    \| \mathcal{G} \|_{\infty}:= \sup_{\omega \in \mathbb{R}} \| \mathcal{G}(\iota\omega) \|_2,
\end{equation}
where $\|\cdot \|_2$ is the spectral norm, namely the largest singular value  $\sigma_{\rm max}(\cdot)$; see \cref{rem:sigma}. 

Boyd, Balakrishnan and Kabamba \cite[Theorem 1]{boyd1989bisection} established a connection between the $H_{\infty}$ norm of a transfer  matrix and locations of eigenvalues of a related Hamiltonian matrix; see also  \cite[Lemma 4.7]{ZDG96}. The following result shows that this connection also holds in the tensor case.

\begin{lemma}
	\label{Pro: H_norm}
	Let $\gamma > 0$ and $\mathcal{A}$ have no imaginary $U$-eigenvalues. Then $\|\mathcal{G}\|_{\infty} <\gamma$ if and only if $\sigma_{\rm max}(\mathcal{D})<\gamma$ and the Hamiltonian tensor $\mathcal{M}_{\gamma}$ has no purely imaginary U-eigenvalues, where 
    \begin{equation} \label{eq:M_r}	\mathcal{M}_{\gamma}:=\left[\begin{matrix}\mathcal{A}+\mathcal{B} *\mathcal{R}^{-1}* \mathcal{D}^{\rm H}* C &  \mathcal{B} *\mathcal{R}^{-1}* \mathcal{B}^{\rm H} \\
	- \mathcal{C}^{\rm H}*\left(\mathcal{I}+\mathcal{D}* \mathcal{R}^{\rm H}* \mathcal{D}^{\rm H}\right)* \mathcal{C} & -\left(\mathcal{A}+\mathcal{B}* \mathcal{R}^{-1} *\mathcal{D}^{\rm H}* \mathcal{C}\right)^{\rm H}
		\end{matrix}\right], 
    \end{equation} 
	and $\mathcal{R}=\gamma^2 \mathcal{I}-\mathcal{D}^{\rm H} * \mathcal{D}$.
\end{lemma}
\begin{proof}
	The proof is an extension of that of \cite[Theorem 1]{boyd1989bisection} for the matrix case. 
	Firstly, we show that $\gamma > 0$ is a singular value of $\mathcal{G}(\iota \omega)$  for some $\omega \in \mathbb{R}$ if and only if $\iota \omega$ is a U-eigenvalue of $\mathcal{M}_{\gamma}$. 
	Let $\gamma$ be a singular value of $\mathcal{G}(\iota \omega)$. 
	By \cref{Lem: tensorSVD} of the tensor SVD to $\mathcal{G}(\iota \omega)$, there is a pair of nonzero tensors $\mathcal{U}, \mathcal{V}$ such that
	\begin{align}
		\label{Eq: 13}
		\begin{cases}
			\left( (\sqrt{\gamma}\mathcal{C})*(\iota \omega \mathcal{I} - \mathcal{A})^{-1}*(\frac{1}{\sqrt{\gamma}}\mathcal{B}) + \mathcal{D} \right)*\mathcal{U} = \gamma \mathcal{V},  \\
			\left( (\frac{1}{\sqrt{\gamma}}\mathcal{B})^{\rm H}*(-\iota \omega \mathcal{I} - \mathcal{A}^{\rm H})^{-1}* (\sqrt{\gamma} \mathcal{C})^{\rm H} + \mathcal{D}^{\rm H} \right)* \mathcal{V} = \gamma \mathcal{U}.
		\end{cases}
	\end{align}
	As $\mathcal{A}$ has no imaginary U-eigenvalues,  $\iota \omega \mathcal{I}-\mathcal{A}$ is invertible. Define nonzero tensors $\mathcal{S}_1, \mathcal{S}_2$ which solve
	\begin{align}
		\label{Eq: 14}
		\begin{cases}
			(\iota \omega \mathcal{I}-\mathcal{A})*\mathcal{S}_1 = \frac{1}{\sqrt{\gamma}}\mathcal{B}*\mathcal{U},  \\
			(-\iota \omega \mathcal{I} - \mathcal{A}^{\rm H})*\mathcal{S}_2 = \sqrt{\gamma} \mathcal{C}^{\rm H} * \mathcal{V}, 
		\end{cases}
	\end{align}
	respectively. 
	Then using tensor operations given in \Cref{sec: pre} in particular  block tensor operations given in \cref{Pro: blocktensor},  we can  rewrite \cref{Eq: 13} as
	\begin{align}
		\label{Eq: 15}
		\begin{bmatrix}
			\mathcal{U} \\
			\mathcal{V}
		\end{bmatrix} = \begin{bmatrix}
			-\mathcal{D} & \gamma \mathcal{I} \\
			\gamma \mathcal{I}  & -\mathcal{D}^{\rm H}
		\end{bmatrix}^{-1} * \begin{bmatrix}
			\sqrt{\gamma}\mathcal{C}  & \mathcal{O} \\
			\mathcal{O}  & \frac{1}{\sqrt{\gamma}}\mathcal{B}^{\rm H}
		\end{bmatrix}*\begin{bmatrix}
			\mathcal{S}_1 \\
			\mathcal{S}_2
		\end{bmatrix} .  	
	\end{align}
	From \cref{Eq: 14} and \cref{Eq: 15}, we obtain
    \begin{multline*}
		\biggl( \begin{bmatrix}
			\mathcal{A}  & \mathcal{O} \\
			\mathcal{O}  & -\mathcal{A}^{\rm H}
		\end{bmatrix} + \begin{bmatrix}
			\frac{1}{\sqrt{\gamma}}\mathcal{B} & \mathcal{O} \\
			\mathcal{O} & -\sqrt{\gamma}\mathcal{C}^{\rm H}
		\end{bmatrix}* \begin{bmatrix}
			-\mathcal{D} & \gamma \mathcal{I} \\
			\gamma \mathcal{I} & - \mathcal{D}^{\rm H}
		\end{bmatrix}^{-1} \\ * \begin{bmatrix}	\sqrt{\gamma}\mathcal{C} & \mathcal{O} \\
			\mathcal{O} & \frac{1}{\sqrt{\gamma}}\mathcal{B}^{\rm H}
		\end{bmatrix} \biggr)*\begin{bmatrix}
			\mathcal{S}_1 \\
			\mathcal{S}_2
		\end{bmatrix} = \iota \omega \begin{bmatrix}
			\mathcal{S}_1 \\
			\mathcal{S}_2
		\end{bmatrix},
    \end{multline*}
    which in fact is  $\mathcal{M}_{\gamma} *\left[\begin{smallmatrix}
		\mathcal{S}_1 \\
		\mathcal{S}_2
	\end{smallmatrix}\right]= \iota \omega \left[\begin{smallmatrix}
		\mathcal{S}_1 \\
		\mathcal{S}_2
	\end{smallmatrix}\right]$. 
	Thus, $\iota \omega$ is a U-eigenvalue of $\mathcal{M}_{\gamma}$. The inverse process of the above analysis also holds.
	
	Next, we proved the result. (Necessity). If $\| \mathcal{G} \|_{\infty} < \gamma$, we have $\sigma_{\rm max} (\mathcal{D}) = \lim_{\omega \to \infty} \sigma_{\rm max} \left(\mathcal{G}(\iota \omega) \right) \leq \| \mathcal{G} \|_{\infty} < \gamma$.  
	Moreover, as $\gamma$ is not a  singular value of $\mathcal{G}(\iota \omega)$,  by the proof given above, $\mathcal{M}_{\gamma}$ has no imaginary U-eigenvalues.  (Sufficiency). Suppose $\sigma_{\rm max} (\mathcal{D})  < \gamma$ and the Hamiltonian tensor $\mathcal{M}_{\gamma}$ has no imaginary U-eigenvalues, but $\| \mathcal{G} \|_{\infty} \ge \gamma$. Since $\sigma_{\rm max} (\mathcal{G}(\iota \omega))$ is a continuous function of $\omega$, and $ \lim_{\omega \to \infty} \sigma_{\rm max} \left(\mathcal{G}(\iota \omega) \right)$ $ =\sigma_{\rm max} (\mathcal{D}) < \gamma$,   there exists an $\omega_0\in \mathbb{R}$ such that $\sigma_{\rm max}(\mathcal{G}(\iota \omega_0)) = \gamma$. Consequently, $\iota \omega_0$ is a U-eigenvalue of $\mathcal{M}_{\gamma}$, which contradicts that   $\mathcal{M}_{\gamma}$ has no imaginary U-eigenvalues.
\end{proof}

\subsection{Unique positive semidefinite solution of the ARTE \cref{ARTE}}
In this  subsection, we establish the unique existence of the positive semidefinite solution of the  ARTE \cref{ARTE}.

We extend stabilizability and detectability to the MLTI system \cref{MLTI-continuous}.

\begin{definition}[stabilizability]\label{def:stabilizability}\rm
	The tensor pair $(\mathcal{A}, \mathcal{B})$ is said to be stabilizable, if for any $\lambda \in \{ z\in \mathbb{C} : \mathrm{Re}(z) \ge 0 \}$ and $\mathcal{Y}$ such that
	$ \mathcal{Y}^{\rm H}*\mathcal{A} = \lambda \mathcal{Y}^{\rm H}$ and $\mathcal{Y}^{\rm H}*\mathcal{B} = \mathcal{O}$, then $\mathcal{Y} = \mathcal{O}$.
\end{definition}
\begin{definition}[detectability]\label{def:detectability}\rm
	The tensor pair $(\mathcal{C}, \mathcal{A})$ is said to be detectable if $(\mathcal{A}^{\rm H}, \mathcal{C}^{\rm H})$ is stabilizable, i.e., if $\lambda \in \{ z\in \mathbb{C} : \mathrm{Re}(z) \ge 0 \}$ and $\mathcal{Y}$ satisfy
	$ \mathcal{A}*\mathcal{Y} = \lambda \mathcal{Y}$ and $\mathcal{C}*\mathcal{Y} = \mathcal{O}$, then $\mathcal{Y} = \mathcal{O}$. 
\end{definition}

\begin{lemma}
	\label{Lem: NoimaginaryEigen}
	Assume that $(\mathcal{A}, \mathcal{B})$ is stabilizable and $(\mathcal{C}, \mathcal{A})$ is detectable. Let $\mathcal{G} = \mathcal{B} * \mathcal{B}^{\mathrm{H}}$ and $\mathcal{K} = \mathcal{C}^{\mathrm{H}} * \mathcal{C}$. Then the Hamiltonian tensor $\mathcal{M}= \left[\begin{smallmatrix}
		\mathcal{A} & \mathcal{G} \\
		\mathcal{K} & -\mathcal{A}^{\mathrm{H}}
	\end{smallmatrix}\right]$ has no purely imaginary U-eigenvalues.
\end{lemma}
\begin{proof}
	We proved this result by contradiction. Suppose that there is a purely imaginary U-eigenvalue $\lambda \in \mathbb{C}$ such that
	$\mathcal{M}* \left[\begin{smallmatrix}
		\mathcal{X} \\
		\mathcal{Y}
	\end{smallmatrix}\right] = \lambda \left[\begin{smallmatrix}
		\mathcal{X} \\
		\mathcal{Y}
	\end{smallmatrix}\right]$, where $  \left[\begin{smallmatrix}
		\mathcal{X} \\
		\mathcal{Y}
	\end{smallmatrix}\right] \ne \mathcal{O}$. 
	By Property 2 of \cref{prop: Hamiltonian}, we have $\mathcal{X}^{\mathrm{H}}* \mathcal{K}* \mathcal{X} + \mathcal{Y}^{\mathrm{H}}* \mathcal{G}* \mathcal{Y} = 0$. 
	As $\mathcal{G} = \mathcal{B} * \mathcal{B}^{\mathrm{H}}$ and $\mathcal{K} = \mathcal{C}^{\mathrm{H}} * \mathcal{C}$ are positive semidefinite, the above equation implies that  $\mathcal{C}* \mathcal{X} = \mathcal{O}\ \mathrm{and} \ \mathcal{B}^{\mathrm{H}} * \mathcal{Y} = \mathcal{O}$. 
	Moreover,  as $\mathrm{Re}(\lambda) = 0$, it follows from Property 4 of \cref{prop: Hamiltonian} that $\mathcal{A}* \mathcal{X} = \lambda \mathcal{X},\ \mathrm{and} \  \mathcal{A}^{\mathrm{H}} * \mathcal{Y} = -\lambda \mathcal{Y}$, which contradict either  $(\mathcal{A}, \mathcal{B})$ being stabilizable when $\mathcal{X} \ne \mathcal{O}$ or $(\mathcal{C}, \mathcal{A})$ being detectable when $\mathcal{Y} \ne \mathcal{O}$.
\end{proof}

We are ready to solve the ARTE \cref{ARTE}, where $\mathcal{G} = \mathcal{B} * \mathcal{B}^{\mathrm{H}}$ and $\mathcal{K} = \mathcal{C}^{\mathrm{H}} * \mathcal{C}$. 
\begin{theorem}
	\label{Thm: ARTEunique}
	Suppose that $(\mathcal{A}, \mathcal{B})$ is stabilizable and $(\mathcal{C}, \mathcal{A})$ is detectable. Let $\mathcal{G} = \mathcal{B} * \mathcal{B}^{\mathrm{H}}$ and $\mathcal{K} = \mathcal{C}^{\mathrm{H}} * \mathcal{C}$. Then there exists a unique positive semidefinite tensor $\mathcal{E}\in \mathbb{C}^{I_1 \times I_1 \times \cdots \times I_N \times I_N}$ solving the ARTE \cref{ARTE} such that $\mathcal{A} - \mathcal{G}* \mathcal{E}$ is stable.
\end{theorem}
\begin{proof}
    \Cref{Lem: NoimaginaryEigen} ensures that the  Hamiltonian tensor $\mathcal{M}$ has no purely  imaginary U-eigenvalues. From \cref{thm:Schur_Hamiltonian} and \cref{Pro: Symblock}, there is a unitary symplectic tensor  $\mathcal{Q} = \left[\begin{smallmatrix}
		\mathcal{Q}_{1} & \mathcal{Q}_{2}\\
		-\mathcal{Q}_{2} & \mathcal{Q}_{1}
	\end{smallmatrix}\right] \in \mathbb{C}^{2I_1 \times 2I_1 \times I_2 \times I_2 \times \cdots \times I_N  \times I_N}$	such that the Hamiltonian tensor $\mathcal{M}$ has the Schur-Hamiltonian tensor decomposition 
	\begin{equation}\label{Eq: SchurHamilton}
		\begin{bmatrix}
			\mathcal{A} & \mathcal{G} \\
			\mathcal{K} & -\mathcal{A}^{\mathrm{H}}
		\end{bmatrix} * \begin{bmatrix}
			\mathcal{Q}_{1} & \mathcal{Q}_{2}\\
			-\mathcal{Q}_{2} & \mathcal{Q}_{1}
		\end{bmatrix} =\begin{bmatrix}
			\mathcal{Q}_{1} & \mathcal{Q}_{2}\\
			-\mathcal{Q}_{2} & \mathcal{Q}_{1}
		\end{bmatrix}*\begin{bmatrix}
			\mathcal{T}  & \mathcal{R} \\
			\mathcal{O}  & -\mathcal{T}^{\mathrm{H}}
		\end{bmatrix}, 
	\end{equation}
	where $\mathcal{T}$ is stable. It follows from \cref{Pro: blocktensor} that
	\begin{align}\label{13}
		\left\{\begin{array}{l}
			\mathcal{A}* \mathcal{Q}_1 - \mathcal{G}* \mathcal{Q}_2 = \mathcal{Q}_1 * \mathcal{T}, \\
			\mathcal{K}* \mathcal{Q}_1 + \mathcal{A}^{\mathrm{H}}* \mathcal{Q}_2 = -\mathcal{Q}_2 * \mathcal{T}.
		\end{array}\right.
	\end{align}
	If the tensor $\mathcal{Q}_1$ is invertible, then eliminating $\mathcal{T}$ from   \cref{13}  yields that $\mathcal{E} := \mathcal{Q}_2* \mathcal{Q}_1^{-1}$ is a solution of the  ARTE \cref{ARTE}, and  $ \mathcal{A} - \mathcal{G}*\mathcal{E} = \mathcal{Q}_1 * \mathcal{T} * \mathcal{Q}_1^{-1}$, which is stable since $\mathcal{T}$ is stable.  Furthermore, by \cref{Pro: Symblock}, $\mathcal{Q}_1^{\rm H} *\mathcal{Q}_2 = \mathcal{Q}_2^{\rm H}* \mathcal{Q}_1$, equivalently $\mathcal{Q}_2* \mathcal{Q}_1^{-1} = (\mathcal{Q}_1^{-1})^{\rm H} * \mathcal{Q}_2^{\rm H}$, which means that the solution $\mathcal{E} = \mathcal{Q}_2* \mathcal{Q}_1^{-1}$ is a Hermitian tensor. Next, we show that $\mathcal{Q}_1$ is invertible. We complete the proof in a similar way as \cite[Theorem 4.1]{paige1981schur}. 
	According to \cref{Thm: SymSVD}, there are unitary tensors $\mathcal{U}$ and $\mathcal{V}$ such that
	\begin{align}\label{10}
		\begin{bmatrix}
			\mathcal{Q}_1  & \mathcal{Q}_2 \\
			-\mathcal{Q}_2  & \mathcal{Q}_1
		\end{bmatrix}  = \begin{bmatrix}
			\mathcal{U}  & \mathcal{O} \\
			\mathcal{O} & \mathcal{U}
		\end{bmatrix}*\begin{bmatrix}
			\mathcal{S} & \mathcal{D} \\
			-\mathcal{D} & \mathcal{S}
		\end{bmatrix}* \begin{bmatrix}
			\mathcal{V}^{\mathrm{H}}  & \mathcal{O}\\
			\mathcal{O} & \mathcal{V}^{\mathrm{H}}
		\end{bmatrix}, 
	\end{align}
	where $\mathcal{S}$ and $\mathcal{D}$ are diagonal tensors, satisfying  $\mathcal{S}_{i_1i_1\cdots i_ni_N} \ge 0$ and $(\mathcal{S}_{i_1i_1\cdots i_Ni_N})^2 + (\mathcal{D}_{i_1i_1\cdots i_Ni_N})^2 = 1$. Moreover, by equation \cref{10} we know that $\mathcal{Q}_{1} = \mathcal{U}* \mathcal{S} * \mathcal{V}^{\mathrm{H}}$ and $ \mathcal{Q}_{2} = \mathcal{U}* \mathcal{D}*\mathcal{V}^{\mathrm{H}}$. Let $\hat{\mathcal{A}} = \mathcal{U}^{\mathrm{H}} * \mathcal{A} * \mathcal{U}$, $\hat{\mathcal{G}} = \mathcal{U}^{\mathrm{H}}* \mathcal{G} * \mathcal{U}$, $\hat{\mathcal{K}} = \mathcal{U}^{\mathrm{H}}* \mathcal{K} * \mathcal{U}$. According to  \cref{Eq: SchurHamilton} and \cref{10}, we have	
	\begin{align}\label{11}
		\begin{bmatrix}
			\hat{\mathcal{A}} & \hat{\mathcal{G}} \\
			\hat{\mathcal{K}} &  -\hat{\mathcal{A}}^{\mathrm{H}}
		\end{bmatrix} * \begin{bmatrix}
			\mathcal{S} & \mathcal{D} \\
			-\mathcal{D} & \mathit{S}
		\end{bmatrix} = \begin{bmatrix}
			\mathcal{S} & \mathcal{D} \\
			-\mathcal{D} & \mathcal{S}
		\end{bmatrix} *\begin{bmatrix}
			\hat{\mathcal{T}}  & \hat{\mathcal{R}} \\
			\mathcal{O}  & -\hat{\mathcal{T}}^{\mathrm{H}}
		\end{bmatrix},
	\end{align}
    where $\hat{\mathcal{R}} = \mathcal{V}^{\mathrm{H}} * \mathcal{R} * \mathcal{V}$, and the tensor $\hat{\mathcal{T}} = \mathcal{V}^{\mathrm{H}} * \mathcal{T} * \mathcal{V}$ have the same U-eigenvalues as $\mathcal{T}$. Let $\hat{A} = \phi(\hat{\mathcal{A}}), \hat{G} = \phi(\hat{\mathcal{G}}), \hat{K} = \phi(\hat{\mathcal{A}}), S = \phi(\mathcal{S})$, and $D = \phi(\mathcal{D})$. From \cref{Lemma: unfolding}, performing $\phi$ on both sides of \cref{11} gives
	\begin{align}\label{12}
		\begin{bmatrix}
			\hat{A} & \hat{G} \\
			\hat{K} &  -\hat{A}^{\mathrm{H}}
		\end{bmatrix} \begin{bmatrix}
			S & D \\
			-D & S
		\end{bmatrix} = \begin{bmatrix}
			S & D \\
			-D & S
		\end{bmatrix}\begin{bmatrix}
			\hat{T}  & \hat{R} \\
			O  & -\hat{T}^{\mathrm{H}}
		\end{bmatrix},
	\end{align}
	where $\hat{T} = \phi(\hat{\mathcal{T}}), \hat{R} = \phi(\hat{\mathcal{R}})$. If the diagonal matrix $S$ is singular, according to the proof of \cite[Theorem 4.1]{paige1981schur}, there exist $\mathbf{w} \ne \mathbf{0} \in \mathbb{C}^{|\mathbf{I}|} $ and $\lambda $ with $\mathrm{Re}(\lambda) < 0$ such that $\hat{A}^{\rm H} \mathbf{w} = - \lambda \mathbf{w}$ and $\hat{G}\mathbf{w} = \mathbf{0}$. Let $\mathcal{W} \in \mathbb{C}^{I_1\times \cdots \times I_N}$ with $\mathcal{W}_{i_1\cdots i_N} = \mathbf{w}_{\mathrm{ivec(\mathbf{i}, \mathbf{I})}}$. Then from the above equation we have	$ \mathcal{A}^{\mathrm{H}}* (\mathcal{U}*\hat{\mathcal{W}}) = -\lambda \mathcal{U}*\hat{\mathcal{W}},\ \mathcal{B}* (\mathcal{U}* \hat{\mathcal{W}}) =\mathcal{O}$ 
	with $\mathrm{Re}(-\lambda) > 0$, which contradicts the stabilizability of $(\mathcal{A}, \mathcal{B})$.
	Therefore, $S$ is nonsingular.  Consequently, the tensor $\mathcal{Q}_1 = \mathcal{U}*\phi^{-1}(S)*\mathcal{V}^{\mathrm{H}}$ is  invertible. 
 
  Finally, rewrite the ARTE \cref{ARTE} as
 \[
    (\mathcal{A}-\mathcal{G}* \mathcal{E})^{\mathrm{H}}* \mathcal{E} + \mathcal{E}* (\mathcal{A}-\mathcal{G}* \mathcal{E}) + \mathcal{E}*\mathcal{G}* \mathcal{E} + \mathcal{K} = \mathcal{O}.
 \]
 As $\mathcal{A} - \mathcal{G}* \mathcal{E}$ is stable, by \cref{lem: LTE}, $\mathcal{E}$ is positive semidefinite and  also unique.
\end{proof}

\Cref{Thm: ARTEunique} is the MLTI version of \cite[Corollary 13.8]{ZDG96}. Based on \cref{lem: LTE} and \cref{Pro: H_norm}, a tensor version of the {\it bounded real lemma} (see for example \cite[Theorem 3.9]{Kimura96} or \cite[Corollary 13.24]{ZDG96}) can  be established.
\begin{theorem}\label{thm:BRL} 
For the continuous-time MLTI system \cref{MLTI-continuous},  let $\mathcal{G}(s)$ be  transfer tensor in Eq. \cref{eq:tf}, and $\mathcal{M}_{\gamma}$ be that in Eq. \cref{eq:M_r}. Assume that $\mathcal{A} \in \mathbb{C}^{I_1 \times I_1 \times \cdots \times I_N \times I_N}$ is stable,  $(\mathcal{A}, \mathcal{B})$ is stabilizable, and $(\mathcal{C}, \mathcal{A})$ is detectable.
	Then the following conditions are equivalent:
	\begin{description}
	    \item[(i)] $\|\mathcal{G}\|_\infty<\gamma$.
		\item[(ii)] $\sigma_{\rm max}(\mathcal{D})<\gamma$ and the Hamiltonian tensor $\mathcal{M}_{\gamma}$ has no purely imaginary U-eigenvalues.
		\item[(iii)] $\sigma_{\rm max}(\mathcal{D})<\gamma$ and there exists a unique positive semidefinite tensor $\mathcal{E}$ that solves
    \begin{align}\label{eq:ARME_M_r}
		  \mathcal{E}*(\mathcal{A}+\mathcal{B} *\mathcal{R}^{-1}* \mathcal{D}^{\rm H}* \mathcal{C}) + (\mathcal{A}+\mathcal{B} *\mathcal{R}^{-1}* \mathcal{D}^{\rm H}* \mathcal{C})^\mathrm{H}*\mathcal{E} +  \\ \mathcal{E}*\mathcal{B} *\mathcal{R}^{-1}* \mathcal{B}^{\rm H} *\mathcal{E} +\mathcal{C}^{\rm H}*\left(\mathcal{I}+\mathcal{D}* \mathcal{R}^{\rm H}* \mathcal{D}^{\rm H}\right)* \mathcal{C} = \mathcal{O}, \nonumber
		\end{align}
    and $\mathcal{A}+\mathcal{B} *\mathcal{R}^{-1}* (\mathcal{D}^{\rm H}* \mathcal{C} + \mathcal{B}^{\rm H}* \mathcal{E}) $ has no purely imaginary U-eigenvalues.
\end{description}
\end{theorem}
\begin{proof}
The equivalence between (i) and (ii) has been established in \cref{Pro: H_norm}.  In the following, we show the equivalence between (ii) and (iii). 
First, show
(iii) $\Rightarrow$ (ii). By equation \cref{eq:ARME_M_r}, we have
\small 
\begin{eqnarray*}
	&&\left[
	\begin{array}{cc}
		\mathcal{I}     &  \mathcal{O} \\
		-\mathcal{E}      &  \mathcal{I}
	\end{array}
	\right] * \mathcal{M}_r *
	\left[
	\begin{array}{cc}
		\mathcal{I}     &  \mathcal{O} \\
		\mathcal{E}      &  \mathcal{I}
	\end{array}
	\right]   
	\\
	&=&
	\left[\begin{matrix}
		\mathcal{A}+\mathcal{B} *\mathcal{R}^{-1}* (\mathcal{D}^{\rm H}* \mathcal{C} + \mathcal{B}^{\rm H}* \mathcal{E}) &  \mathcal{B} *\mathcal{R}^{-1}* \mathcal{B}^{\rm H} \\
		\mathcal{O} & -\left(\mathcal{A}+\mathcal{B} *\mathcal{R}^{-1}* (\mathcal{D}^{\rm H}* \mathcal{C} + \mathcal{B}^{\rm H}* \mathcal{E})\right)^{\rm H}
	\end{matrix}\right].
\end{eqnarray*}
\normalsize
Thus, (iii) implies (ii). Next, show (ii) $\Rightarrow$ (iii). 
Since $(\mathcal{A}, \mathcal{B})$ is stabilizable,  $(\mathcal{A}+\mathcal{B} *\mathcal{R}^{-1}* \mathcal{D}^{\rm H}* \mathcal{C}, \mathcal{B} *\mathcal{R}^{-1}* \mathcal{B}^{\rm H})$ is also stabilizable.  Then using a procedure similar to that in the proof of \cite[Theorem 13.6]{ZDG96} and then \cite[Theorem 13.5]{ZDG96}, it can be readily shown  that  \cref{eq:ARME_M_r} has a Hermitian solution $\mathcal{E}$ and $ \mathcal{A}+\mathcal{B} *\mathcal{R}^{-1}* (\mathcal{D}^{\rm H}* \mathcal{C} + \mathcal{B}^{\rm H}* \mathcal{E})$ is stable, hence has no purely imaginary U-eigenvalues. Finally, as (ii) is equivalent to (i), but the latter implies that the tensor  $\mathcal{A}+\mathcal{B} *\mathcal{R}^{-1}* \mathcal{D}^{\rm H}* \mathcal{C}$ is stable. As the result, by \cref{lem: LTE}, the solution $\mathcal{E}$ to the Lyapunov tensor equation \cref{eq:ARME_M_r} is positive semidefinite and unique.
\end{proof}

\begin{remark}
	\label{Remark:ARTE}
Setting $\mathcal{D}=\mathcal{O}$ and $\gamma=1$, \cref{eq:ARME_M_r} reduces to the ARTE \cref{ARTE} with $\mathcal{G}=-\mathcal{B}*\mathcal{B}^{\rm H}$.
\end{remark}

\subsection{A tensor version of the small gain theorem} \label{subsec:small gain theorem}

We begin with a lemma.
\begin{lemma} \label{lem:small gain}
Let $\mathcal{G}, \mathcal{H}\in \mathcal{FH}_{\infty} $. Then $(\mathcal{I}-\mathcal{G}*\mathcal{H})^{-1} \in \mathcal{FH}_{\infty}$ if $\|\mathcal{G} \|_\infty \ \|\mathcal{H} \|_\infty <1$.
\end{lemma}
\begin{proof}
Recall the tensor unfolding operation $\phi$ in   \cref{def:phi}.   Actually,  $\phi$ is a group isomorphism.  The case of even-order tensors with $N=2$ was established in  \cite[Lemma 3.1]{Brazell2013}, which can be easily extended to general cases. Therefore,    
\begin{equation} \label{eq:jan27_GH}
  \phi\left(\mathcal{G}(s)*\mathcal{H}(s)\right) = \phi(\mathcal{G}(s))\ \phi(\mathcal{H}(s)).  
\end{equation}
Moreover, recall  the unfolding determinant $\mathrm{det}_U(\cdot)$ defined in the properties below  \cref{def:phi}, we have 
\begin{equation}\label{eq:jan27_gege}
\mathrm{det}_U[(\mathcal{I}-\mathcal{G}(s)*\mathcal{H}(s))^{-1}]  = \mathrm{det}\{\phi[(\mathcal{I}-\mathcal{G}(s)*\mathcal{H}(s))^{-1}]\}
=  \frac1{ \mathrm{det}[\phi(\mathcal{I}-\mathcal{G}(s)*\mathcal{H}(s))]}.
\end{equation}
Let $\mathcal{G}, \mathcal{H}\in \mathcal{FH}_{\infty} $. Clearly,  $\mathcal{G}*\mathcal{H} \in \mathcal{FH}_{\infty} $ too. Due to \cref{eq:jan27_gege}, to show $(\mathcal{I}-\mathcal{G}*\mathcal{H})^{-1} \in \mathcal{FH}_{\infty}$, it suffices to show that $\mathrm{det}(I-\phi(\mathcal{G}(s)* \mathcal{H}(s)))$ has no zeros in the closed right-half plane. However, by \cref{eq:jan27_GH},
\begin{equation} \label{eq:jan27_temp1}
 \|\phi\left(\mathcal{G}(s)*\mathcal{H}(s)\right)\|_2   \leq
  \|\phi(\mathcal{G}(s))\|_2   \ \|\phi(\mathcal{H}(s))\|_2  = \|\mathcal{G}(s)\|_2  \  \|\mathcal{H}(s)\|_2,
\end{equation}
where   $\|\cdot \|_2$ is the spectral norm defined in  \cref{def:spectrum}. If $\|\mathcal{G} \|_\infty \  \|\mathcal{H} \|_\infty <1$, then from \cref{eq:jan27_temp1} we have for each $s$ in the closed right-half plane,
\begin{equation}
 \|\phi\left(\mathcal{G}(s)*\mathcal{H}(s)\right)\|_2 \leq  \|\mathcal{G}(s)\|_2  \  \|\mathcal{H}(s)\|_2 \leq     \|\mathcal{G}\|_\infty\  \|\mathcal{H}\|_\infty <1.
\end{equation}
Consequently, $\mathrm{det}(I-\phi(\mathcal{G}(s)* \mathcal{H}(s)))$ has no zeros in the closed right-half plane.
\end{proof}

\begin{figure}[htbp!]
		\centering
\includegraphics[width=0.45\textwidth]{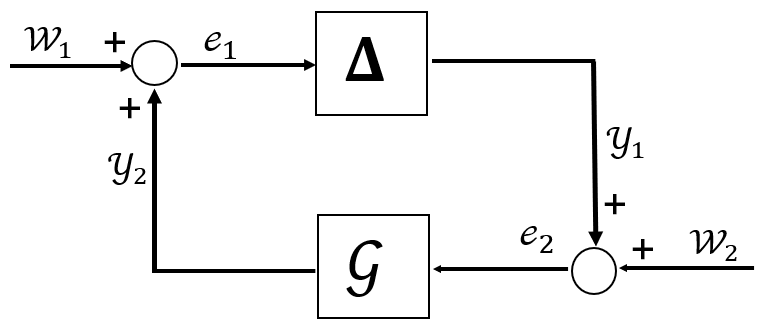}
		\caption{Closed-loop stability analysis, where the state-space realizations   $(\mathcal{A}_\mathcal{G},\mathcal{B}_\mathcal{G},\mathcal{C}_\mathcal{G},\mathcal{D}_\mathcal{G})$ of the plant  $\mathcal{G}$ and $(\mathcal{A}_\Delta,\mathcal{B}_\Delta,\mathcal{C}_\Delta,\mathcal{D}_\Delta)$ of the uncertainty $\Delta $ are stabilizable and detectable, and the corresponding transfer tensors $\mathcal{G}, \Delta \in \mathcal{FH}_{\infty}$.}
 \label{fig_small_gain} 
	\end{figure}

Next, we present the small gain theorem for the feedback system in \cref{fig_small_gain}.

First, we introduce the concept of internal stability of the feedback system in \cref{fig_small_gain}, which is an extension of \cite[Def. 4.4]{SP05}.

\begin{definition}
    The feedback system in \cref{fig_small_gain} is internally stable if neither $\mathcal{G}$ nor $\Delta$ contains  hidden unstable modes and the injection of bounded  external  signals at any place in the feedback system results in bounded output signals measured anywhere in the system.  
\end{definition}

\begin{remark}
 In linear systems  theory,  an LTI system contains no hidden unstable modes if and only if the system is stabilizable and detectable.  In this case, the input-output (BIBO) stability of the transfer matrix is equivalent to  that the $A$-matrix of  its state-space realization is a Hurwitz matrix, i.e., the system is   internally stable.    As given in \cref{rem:stability},  $\lambda$ is  a $U$-eigenvalue of a tensor $\mathcal{A}$ of the MLTI system \cref{MLTI-continuous} if and only if it is an eigenvalue of the associate matrix $\phi(\mathcal{A})$.  Therefore, under the assumption of stabilizability and detectability (see  \cref{def:stabilizability} and \cref{def:detectability}), it can be easily shown that the BIBO stability of the transfer tensor is equivalent to the internal stability of the MLTI system \cref{MLTI-continuous}.
 \end{remark}

We are ready to present a tensor version of the small gain theorem, where the proof of the  necessity part  is an immediate extension of the matrix case given in \cite[Theorem 8.2]{DP00}.

\begin{theorem} \label{thm:small gain}
    Let $\mathcal{G} \in \mathcal{FH}_{\infty}$. Then the closed-loop system in \cref{fig_small_gain}  is internally stable for all $\Delta \in \mathcal{FH}_{\infty}$ with   $\|\Delta \|_\infty \leq 1$ if and only if $\|\mathcal{G}\|_\infty <1$.
\end{theorem}

\begin{proof}  
 ({\it Sufficiency.})  As  $\mathcal{G}, \Delta \in \mathcal{FH}_{\infty}$ and they are stabilizable and detectable, similar to the matrix case studied in \cite[Chapter 5]{ZDG96} in particular \cite[Corollary 5.6]{ZDG96}, the closed-loop system in \cref{fig_small_gain}  is internally stable if and only if $(\mathcal{I}-\mathcal{G}*\Delta)^{-1} \in \mathcal{FH}_{\infty}$. Thus, the sufficient condition follows directly from \cref{lem:small gain}. 
({\it Necessity.})  By contradiction. Assume $\|\mathcal{G}\|_\infty\geq 1$. Then there exists an $\omega_0 \in [0, \infty]$ such that $\|\mathcal{G}(\imath\omega_0)\|_2\geq 1$.  (Notice $\mathcal{G}(\imath \infty )= \mathcal{D}$.)  Define  $\lambda = \|\mathcal{G}(\imath\omega_0)\|_2^2$. Then define a constant tensor $\Delta = \frac{\mathcal{G}(\imath\omega_0)^\mathrm{H}}{\lambda}$.  Clearly, $\|\Delta\|_\infty = \|\Delta\|_2 = \frac1{\|\mathcal{G}(\imath\omega_0)\|_2}\leq 1$. Hence, $\Delta \in \mathcal{FH}_{\infty}$.  On the other hand, the tensor
\[ \mathcal{I}-\mathcal{G}(\imath\omega_0)\ast \Delta = \frac1{\lambda} (\|\mathcal{G}\left(\imath\omega_0)\|_2^2\;  \mathcal{I} -  \mathcal{G}(\imath\omega_0) \ast \mathcal{G}(\imath\omega_0)^\mathrm{H}\right)
\]
is {\it singular}. Consequently, the closed-loop system is not internally stable.
\end{proof}

\subsection{Perturbation analysis} \label{subsec:pertub}

In practical applications, the coefficient tensors $\mathcal{A}, \mathcal{G}, \mathcal{K}$ in the ARTE \cref{ARTE} may have uncertainties $\Delta \mathcal{A}, \Delta \mathcal{G}, \Delta \mathcal{K}$, respectively. We in fact are given a perturbed ARTE
\begin{equation}\label{PerturbedARTE}
	(\mathcal{A}+\Delta \mathcal{A})^{\mathrm{H}}* \tilde {\mathcal{E}}+ \tilde {\mathcal{E}}* (\mathcal{A}+\Delta \mathcal{A}) - \tilde{\mathcal{E}}* (\mathcal{G} + \Delta \mathcal{G}) * \tilde{\mathcal{E}} + \mathcal{K} + \Delta \mathcal{K} = \mathcal{O}.
\end{equation}
It is clearly meaningful and important to measure the perturbation error $\Delta \mathcal{E} = \tilde{\mathcal{E}} - \mathcal{E} $ between the positive semidefinite solutions of the original ARTE \cref{ARTE} and the perturbed one  \cref{PerturbedARTE}. In this subsection, a perturbation analysis under the assumptions in \cref{Thm: ARTEunique} is conducted.

Subtracting \cref{PerturbedARTE} from  \cref{ARTE} and omitting the second-order terms yields
\begin{align} \label{eq:pert_22Jan}
	&(\mathcal{A} - \mathcal{G}*\mathcal{E})^{\rm H}* \Delta \mathcal{E} + \Delta \mathcal{E}*(\mathcal{A} - \mathcal{G}* \mathcal{E}) \\
	\approx &-\Delta \mathcal{A}^{\mathrm{H}} * \mathcal{E} - \mathcal{E}* \Delta \mathcal{A} + \mathcal{E}* \Delta \mathcal{G} * \mathcal{E} - \Delta \mathcal{K}.  \nonumber
\end{align}
 By \cref{Pro: KronVec}, from \cref{eq:pert_22Jan} we obtain that
\begin{align} \label{Eq: Perturb}
	\mathcal{Z}* \mathrm{Vec}(\Delta \mathcal{E}) \approx -(\mathcal{E}^{\top}\otimes \mathcal{I}_{\mathbf{I}}) * \mathrm{Vec}(\Delta \mathcal{A}^{\rm H}) - ( \mathcal{I}_{\mathbf{I}}\otimes \mathcal{E} )* \mathrm{Vec}(\Delta \mathcal{A}) \\ 
	+  (\mathcal{E}^{\top}\otimes \mathcal{E})*\mathrm{Vec}(\Delta \mathcal{G})  - \mathrm{Vec}(\Delta \mathcal{K}), \nonumber
\end{align}
where $\mathcal{Z} = \mathcal{I}_{\mathbf{I}}\otimes (\mathcal{A} - \mathcal{E}* \mathcal{G})^{\mathrm{H}} + (\mathcal{A} - \mathcal{G}* \mathcal{E})^{\top} \otimes \mathcal{I}_{\mathbf{I}}$. According to \cref{Lem: KronUeigen}, $\mathcal{Z}$ is an invertible tensor since  the tensor $\mathcal{A} - \mathcal{G}* \mathcal{E}$ is stable. Left-multiplying the inverse of $\mathcal{Z}$ on both sides of  \cref{Eq: Perturb}  and taking the Frobenius norm, we immediately obtain the first-order perturbation bound.

\begin{theorem}
	\label{Thm: Perturb}
	Under the assumptions in \cref{Thm: ARTEunique}, let $\mathcal{E}$ be the unique  positive semidefinite solution of the ARTE \cref{ARTE}, and $\mathcal{E}+\Delta \mathcal{E}$ be the solution of the perturbed one \cref{PerturbedARTE}. Then
	\begin{align*}
		\| \Delta \mathcal{E} \|_F \lesssim & \left( \| \mathcal{Z}^{-1}* ( \mathcal{I}_{\mathbf{I}}\otimes \mathcal{E}) \|_2 +  \|\mathcal{Z}^{-1}*(\mathcal{E}^{\top}\otimes \mathcal{I}_{\mathbf{I}})  \|_2\right) \| \Delta \mathcal{A} \|_F \\ 
		&  + \| \mathcal{Z}^{-1}* (\mathcal{E}^{\top}\otimes \mathcal{E}) \|_2  \| \Delta \mathcal{G} \|_F 
		+ \| \mathcal{Z}^{-1} \|_2 \| \Delta \mathcal{K} \| _F,
	\end{align*}
	where $\lesssim$ means that ``$\le$'' holds up to the first order in $\| \Delta \mathcal{E} \|_F$.
\end{theorem}

Based on the condition number theory  \cite{rice1966theory, zhou2009perturbation, ghavimi1995backward}, we consider three kinds of normwise condition numbers for ARTEs, which are 
\[
\kappa_i=\lim _{\epsilon \rightarrow 0} \sup _{\Delta_i \leq \epsilon} \frac{\|\Delta \mathcal{E}\|_F}{\epsilon\| \mathcal{E} \|_F} \quad (i=1,2,3),
\]
where
\begin{align*}
	\begin{cases}
		\Delta_1 = 
		\left \| \begin{bmatrix}
			\frac{\|\Delta \mathcal{A} \|_F}{\| \mathcal{A} \|_F} & \frac{\|\Delta \mathcal{G} \|_F}{\| \mathcal{G} \|_F} & \frac{\|\Delta \mathcal{K} \|_F}{\| \mathcal{K} \|_F} 
		\end{bmatrix}  \right \|_F , \\
		\Delta_2 =\max \left\{\frac{\|\Delta \mathcal{A}\|_F}{\| \mathcal{A} \|_F}, \frac{\|\Delta \mathcal{G}\|_F}{\| \mathcal{G} \|_F}, \frac{\|\Delta \mathcal{K} \|_F}{\| \mathcal{K} \|_F}\right\}, \\
		\Delta_3 = \frac{\left \| \begin{bmatrix}
				\Delta \mathcal{A} & \Delta \mathcal{G} & \Delta \mathcal{K}
			\end{bmatrix} \right \|_F }{\left \| \begin{bmatrix}
				\mathcal{A} & \mathcal{G} & \mathcal{K}
			\end{bmatrix} \right \|_F},
	\end{cases}
\end{align*}
they can be regarded as three different measures of the relative error of the perturbations. 

\begin{theorem}
	\label{Thm: ConNum}
	Under the assumptions in \cref{Thm: ARTEunique}, let $\mathcal{E}$ be the unique positive semidefinite solution of the ARTE \cref{ARTE}. Then  
	\begin{align*}
		\begin{cases}
			\kappa_1 \lesssim \frac{\left\|\mathcal{Z}^{-1} * \mathcal{S}_1\right\|_2}{\|\mathcal{E}\|_F}, \\
			\kappa_2 \lesssim \min \left\{\sqrt{3} \kappa_1, \frac{\eta_c}{\|\mathcal{E}\|_F}\right\}, \\
			\kappa_3 \lesssim \frac{\left\|\mathcal{Z}^{-1} * \mathcal{S}_2\right\|_2 \sqrt{\|\mathcal{A}\|_F^2+\|\mathcal{G}\|_F^2+\|\mathcal{K}\|_F^2}}{\|\mathcal{E}\|_F},
		\end{cases}
	\end{align*}
    where tensors 
	\begin{align*}
		\begin{cases}
			\mathcal{S}_1 =\begin{bmatrix}
				\alpha \| \mathcal{A} \|_F(\mathcal{E}^{\top}\otimes \mathcal{I}_{\mathbf{I}})  & \beta \| \mathcal{A} \|_F (\mathcal{I}_{\mathbf{I}}\otimes \mathcal{E} )  & \| \mathcal{G} \|_F \mathcal{E}^{\top}\otimes \mathcal{E} & \| \mathcal{K} \|_F\mathcal{I}_{\mathbf{I}\odot \mathbf{I}} 
			\end{bmatrix}, \\
			\mathcal{S}_2 = \begin{bmatrix}
				\alpha (\mathcal{E}^{\top}\otimes \mathcal{I}_{\mathbf{I}})  & \beta  (\mathcal{I}_{\mathbf{I}}\otimes \mathcal{E} )  & \mathcal{E}^{\top}\otimes \mathcal{E} & \mathcal{I}_{\mathbf{I}\odot \mathbf{I}} 
			\end{bmatrix}
		\end{cases}
	\end{align*}
	with $\frac{1}{\alpha^2} + \frac{1}{\beta ^2} = 1$, and
    \begin{multline*}
            \eta_c=\left( \| \mathcal{Z}^{-1}* ( \mathcal{I}_{\mathbf{I}}\otimes \mathcal{E}) \|_2 + \|\mathcal{Z}^{-1}*(\mathcal{E}^{\top}\otimes \mathcal{I}_{\mathbf{I}})  \|_2\right) \| \mathcal{A} \|_F \\ 
		+ \| \mathcal{Z}^{-1}* (\mathcal{E}^{\top}\otimes \mathcal{E}) \|_2  \| \mathcal{G} \|_F + \| \mathcal{Z}^{-1} \|_2 \| \mathcal{K} \|_F. 
    \end{multline*}
\end{theorem}
\begin{proof}
    By block tensor techniques presented in \Cref{sec:tensor_blocking}, we can rewrite \cref{Eq: Perturb} as
        \begin{align*}
            \mathcal{Z} * \mathrm{Vec}(\Delta \mathcal{E}) & \approx \mathcal{S}_2 * \begin{bmatrix}
            -\frac{1}{\alpha} \mathrm{Vec}({\Delta \mathcal{A}}^{\rm H})^{\top} & -\frac{1}{\beta}\mathrm{Vec}(\Delta \mathcal{A})^{\top} & \mathrm{Vec}(\Delta \mathcal{G})^{\top} & -\mathrm{Vec}(\Delta\mathcal{K})^{\top}
        \end{bmatrix}^{\top}.
        \end{align*}
    Left-multiplying the inverse of $\mathcal{Z}$ and taking the Frobenius norm of both sides yields
	\begin{align*} 
		\|  \Delta \mathcal{E} \|_F & \lesssim
		\| \mathcal{Z}^{-1}* \mathcal{S}_2\|_2\sqrt{\|\Delta \mathcal{A}\|_F^2+\|\Delta \mathcal{G}\|_F^2+\|\Delta \mathcal{K}\|_F^2} \\
		& = \| \mathcal{Z}^{-1}* \mathcal{S}_2\|_2\sqrt{\|\mathcal{A}\|_F^2+\|\mathcal{G}\|_F^2+\|\mathcal{K}\|_F^2}\ \Delta_3. 
	\end{align*}
	Hence,
	$ \kappa_3 \lesssim \frac{\left\|\mathcal{Z}^{-1} * \mathcal{S}_2\right\|_2 \sqrt{\|\mathcal{A}\|_F^2+\|\mathcal{G}\|_F^2+\|\mathcal{K}\|_F^2}}{\|\mathcal{E}\|_F}$. 
	Similarly, \cref{Eq: Perturb} can be rewritten as
    \begin{align*}
        \mathcal{Z} * \mathrm{Vec}(\Delta \mathcal{E}) & \approx  \mathcal{S}_1 *\begin{bmatrix}
            -\frac{\mathrm{Vec}({\Delta \mathcal{A}}^{\rm H}) ^{\top}}{\alpha \| \mathcal{A} \|_F} & -\frac{\mathrm{Vec}(\Delta \mathcal{A})^{\top}}{\beta \| \mathcal{A} \|_F} & \frac{\mathrm{Vec}(\Delta \mathcal{G})^{\top}}{\| \mathcal{G} \|_F} & -\frac{\mathrm{Vec}(\Delta\mathcal{K})^{\top}}{\|\mathcal{K}\|_F}
        \end{bmatrix}^{\top}.
    \end{align*}
It is easily verified that 
$\|\Delta \mathcal{E}\|_F  \lesssim\left\| \mathcal{Z}^{-1}* \mathcal{S}_1\right\|_2 \sqrt{\frac{\|\Delta \mathcal{A} \|_F^2}{\| \mathcal{A} \|_F^2}+\frac{\|\Delta \mathcal{G}\|_F^2}{\| \mathcal{G} \|_F^2}+\frac{\|\Delta \mathcal{K}\|_F^2}{\| \mathcal{K} \|_F^2}}.
$ Thus,
	$  \kappa_1 \lesssim \frac{\left\|\mathcal{Z}^{-1} * \mathcal{S}_1\right\|_2}{\|\mathcal{E}\|_F}$.  
	From the definitions of the relative errors $\Delta_1$ and $\Delta_2$, we have  $\Delta_1 \le \sqrt{3} \Delta_2$. Besides, it follows from \cref{Thm: Perturb} that $	\|\Delta \mathcal{E}\|_F \lesssim \eta_c \Delta_2$.
	As a result, we get the upper bounds for $\kappa _2$.
\end{proof}

If the perturbations $\Delta \mathcal{A}$ on the coefficient tensor $\mathcal{A}$ are real, then tighter perturbation bounds can be given. Denote by $\mathcal{P}^{\mathbf{i}\wedge \mathbf{j}} \in \mathbb{C}^{I_1 \times I_1 \times \cdots \times I_N \times I_N}$ the tensor that has only one nonzero entry $(\mathcal{P}^{\mathbf{i}\wedge \mathbf{j}})_{i_1j_1\cdots i_Nj_N} = 1$, and define the permutation tensor $\mathcal{P} = \sum_{\mathbf{1}\le \mathbf{i}\le \mathbf{I}} \sum_{\mathbf{1}\le \mathbf{j}\le \mathbf{I}} \left(\mathcal{P}^{\mathbf{i}\wedge \mathbf{j}}\right)^{\top} \otimes \mathcal{P}^{\mathbf{i}\wedge \mathbf{j}} \in \mathbb{C}^{I_1^2\times I_1^2\times \cdots \times I_N^2\times I_N^2} $.

\begin{proposition}
	\label{Pro: transposePermutaion}
	Let $\mathcal{E}\in \mathbb{C}^{I_1\times I_1 \times \cdots \times I_N\times I_N}$. Then
	 $\mathrm{Vec}(\mathcal{E}^{\top}) = \mathcal{P}* \mathrm{Vec}(\mathcal{E})$.
\end{proposition}
\begin{proof}
	From the definition of the Einstein product,
	$ \mathcal{E}^{\top} = \sum_{\mathbf{1}\le \mathbf{i}\le \mathbf{I}} \sum_{\mathbf{1}\le \mathbf{j}\le \mathbf{I}}\mathcal{P}^{\mathbf{i}\wedge \mathbf{j}} * \mathcal{E} * \mathcal{P}^{\mathbf{i}\wedge \mathbf{j}}$.  
	From the linearity of Vec operation and \cref{Pro: KronVec}, performing the Vec operation on both sides gives that
	\begin{align*}
		\mathrm{Vec}(\mathcal{E}^{\top}) & = \sum_{\mathbf{1}\le \mathbf{i}\le \mathbf{I}} \sum_{\mathbf{1}\le \mathbf{j}\le \mathbf{I}} \left( (\mathcal{P}^{\mathbf{i}\wedge \mathbf{j}})^{\top} \otimes \mathcal{P}^{\mathbf{i}\wedge \mathbf{j}}\right) * \mathrm{Vec}(\mathcal{E}) = \mathcal{P}* \mathrm{Vec}(\mathcal{E}).
	\end{align*}
	The last equality is obtained from the associative law of the Einstein product.
\end{proof}

According to \cref{Pro: transposePermutaion}, if the perturbation tensor $\Delta \mathcal{A}$ is a real tensor, \cref{Eq: Perturb} reduces to
\begin{multline}
    \mathcal{Z}* \mathrm{Vec}(\Delta \mathcal{E}) \approx - \left( \mathcal{I}_{\mathbf{I}}\otimes \mathcal{E} + (\mathcal{E}^{\top}\otimes \mathcal{I}_{\mathbf{I}})* \mathcal{P}  \right)* \mathrm{Vec}(\Delta \mathcal{A}) \\
	+  (\mathcal{E}^{\top}\otimes \mathcal{E})*\mathrm{Vec}(\Delta \mathcal{G})  - \mathrm{Vec}(\Delta \mathcal{K}). \nonumber
\end{multline}
Thus, we can obtain tighter bounds for the perturbations if $\Delta \mathcal{A}$ is real, which is given in the following remark. 
\begin{remark}
	\label{Remark 4.2}
	When the perturbation $\Delta \mathcal{A}$ is a real tensor, the perturbed error bound can be tighter as given below
	\begin{multline*}
		\| \Delta \mathcal{E} \|_F \lesssim  \| \mathcal{Z}^{-1}* \left( \mathcal{I}_{\mathbf{I}}\otimes \mathcal{E} + (\mathcal{E}^{\top}\otimes \mathcal{I}_{\mathbf{I}})* \mathcal{P}  \right)  \|_2 \| \Delta \mathcal{A} \|_F \\ 
		 + \| \mathcal{Z}^{-1}* (\mathcal{E}^{\top}\otimes \mathcal{E}) \|_2  \| \Delta \mathcal{G} \|_F 
		+ \| \mathcal{Z}^{-1} \|_2 \| \Delta \mathcal{K} \| _F,
	\end{multline*}
	and the tensors $\mathcal{S}_1, \mathcal{S}_2$ and the scalar  $\eta_c$ in \cref{Thm: ConNum} become
	\begin{align*}
		\begin{cases}
			\mathcal{S}_1 =\begin{bmatrix}
				\| \mathcal{A} \|_F \left(\mathcal{I}_{\mathbf{I}}\otimes \mathcal{E} + (\mathcal{E}^{\top}\otimes \mathcal{I}_{\mathbf{I}})* \mathcal{P}\right)  & \| \mathcal{G} \|_F\mathcal{E}^{\top}\otimes \mathcal{E} & \| \mathcal{K} \|_F\mathcal{I}_{\mathbf{I}\odot \mathbf{I}} 
			\end{bmatrix}, \\
			\mathcal{S}_2 = \begin{bmatrix}
				\mathcal{I}_{\mathbf{I}}\otimes \mathcal{E} + (\mathcal{E}^{\top}\otimes \mathcal{I}_{\mathbf{I}})* \mathcal{P}  & \mathcal{E}^{\top}\otimes \mathcal{E} & \mathcal{I}_{\mathbf{I}\odot \mathbf{I}} 
			\end{bmatrix}, 
		\end{cases}
	\end{align*}
	and 
	\begin{multline*}
        \eta_c=\left\|\mathcal{Z}^{-1} *\left(\mathcal{I}_{\mathbf{I}} \otimes \mathcal{E}+(\mathcal{E}^{\top} \otimes \mathcal{I}_{\mathbf{I}}) * \mathcal{P} \right)\right\|_2 \|\mathcal{A} \|_F \\ + \left\|\mathcal{Z}^{-1} *(\mathcal{E}^{\top}\otimes \mathcal{E})\right\|_2 \|\mathcal{G} \|_F + \left\|\mathcal{Z}^{-1}\right\|_2 \|\mathcal{K} \|_F.
	\end{multline*}
\end{remark}

\section{Numerical algorithms and examples}
\label{sec. 5}
In this section, we derive tensor-based numerical algorithms for the ARTE \cref{ARTE} that make use of the multidimensional data structure, and borrow a numerical example in \cite{chen2021multilinear} to demonstrate our results.  

It is well-known that the Newton method is a classical iterative method for solving nonlinear equation $f(x) = 0$, which has the recursive format \cite{kelley2003solving}
\begin{equation}\label{16}
	\begin{cases}
		f'(x_k)\Delta _k  = -f(x_k), \\
		x_{k+1} = x_k + \Delta_k,
	\end{cases}
\end{equation}
where $f'(x)$ is the Fr$\rm \acute{e}$chet derivative \cite{kesavan2022nonlinear}. Let
\begin{equation} \label{eq:f_26Jan}
	f(\mathcal{E}) = \mathcal{A}^{\mathrm{H}}* \mathcal{E} + \mathcal{E}* \mathcal{A} - \mathcal{E}*\mathcal{G}* \mathcal{E} + \mathcal{K}.
\end{equation}
Then under the mild assumptions that $\mathcal{G}$ and $\mathcal{E}$ are Hermitian tensors,
\begin{align*}
	f(\mathcal{E}+\Delta \mathcal{E}) &= \mathcal{A}^{\mathrm{H}} * (\mathcal{E} + \Delta \mathcal{E}) + (\mathcal{E} + \Delta \mathcal{E})* \mathcal{A} - (\mathcal{E} + \Delta \mathcal{E})* \mathcal{G} * (\mathcal{E} + \Delta \mathcal{E}) + \mathcal{K} \\
	& = f(\mathcal{E}) + \mathcal{L} (\Delta \mathcal{E}) - \Delta \mathcal{E} * \mathcal{G}* \Delta \mathcal{E},
\end{align*}
where
$
\mathcal{L}(\Delta \mathcal{E})
= (\mathcal{A} - \mathcal{G}* \mathcal{E})^{\mathrm{H}}* \Delta\mathcal{E} + \Delta \mathcal{E} * (\mathcal{A} - \mathcal{G}* \mathcal{E})$ 
contains all the first order terms. Hence, $f$ is Fr$\rm \acute{e}$chet differentiable and 
\begin{equation} \label{eq:frechet}
	f'(\mathcal{E})(\Delta \mathcal{E}) = \mathcal{L} (\Delta \mathcal{E}).
\end{equation}
Substituting \cref{eq:frechet} into \cref{16} yields the Lyapunov tensor equation \cref{17} given in \cref{alg:1} below. 

Based on the analysis presented above, we obtain the Newton method for solving the ARTE \cref{ARTE} as given in \cref{alg:1}, where in each iteration a Lyapunov tensor equation (LTE) is solved based on the \textbf{Newton-LTE-BiCG} (biconjugate gradient) subroutine or the \textbf{Newton-VecLTE-BiCG} subroutine.

\begin{algorithm}[!htpb]	\renewcommand{\algorithmicrequire}{\textbf{Input:}}
	\renewcommand{\algorithmicensure}{\textbf{Output:}}
	\caption{Tensor-based Newton method for the ARTE \cref{ARTE}}
	\label{alg:1}
	\centering
	\begin{algorithmic}[1]
		\REQUIRE Coefficient tensors $ \mathcal{A}, \mathcal{G}, \mathcal{K} \in \mathbb{C}^{I_1\times I_1 \times \cdots \times I_N \times I_N}$, initial value $\mathcal{E}_0 \in \mathbb{C}^{I_1\times I_1 \times \cdots \times I_N \times I_N}$ and $k = 0$, residual error tolerance $\epsilon$.
		\ENSURE Approximation of the Hermitian  positive semidefinite solution $\mathcal{E} \in \mathbb{C}^{I_1 \times I_1 \times  \cdots \times I_N \times I_N}$.
		
		\STATE Compute $\mathcal{A}_k = \mathcal{A} - \mathcal{G} * \mathcal{E}_k$ and  $\mathcal{K}_k = \mathcal{E}_k^{\rm H}* \mathcal{G} * \mathcal{E}_k + \mathcal{K}$. 
		\STATE Solve $\mathcal{E}_{k+1}$ from the Lyapunov tensor equation
		\begin{equation} \label{17}
			\mathcal{A}_k^{\mathrm{H}} * \mathcal{E}_{k+1} + \mathcal{E}_{k+1}* \mathcal{A}_k + \mathcal{K}_k = \mathcal{O}
		\end{equation}
		by either the Newton-LTE-BiCG or the Newton-VecLTE-BiCG method.
		\STATE Compute the norm residual error $r = \| 	f( \mathcal{E}_{k+1}) \|_F$, if $r < \epsilon $, stop the loops and set $\mathcal{E}_{k+1}$ as the approximation; Otherwise, set $k = k+1$, turn back to step 1. 
	\end{algorithmic}  
\end{algorithm} 

\textbf{Newton-LTE-BiCG}: Adopt the tensor-based method `Algorithm BiCG1' \cite{hajarian2020conjugate} to solve the Lyapunov tensor equation \cref{17} in \cref{alg:1}.

\textbf{Newton-VecLTE-BiCG}: According to \cref{Pro: KronVec}, we rewrite the Lyapunov tensor equation \cref{17} as the equivalent multilinear system $\left( \mathcal{I} \otimes \mathcal{A}_k^{\mathrm{H}} + \mathcal{A}_k^{\top} \otimes \mathcal{I} \right)* \mathcal{Y}_{k+1} = - \mathrm{Vec}(\mathcal{K}_k)$.
This is VecLTE part. Then solve $\mathcal{Y}_{k+1} $ by the tensor BiCG method \cite{huang2021numerical}, and set $\mathcal{E}_{k+1} = \mathrm{Vec}^{-1}(\mathcal{Y}_{k+1})$.

Nowadays, tensor products can be performed efficiently on modern computational platforms without matricization. For example, in  \cite{huang2021numerical} the authors demonstrated the advantage of solving the multilinear system  $\mathcal{A}* \mathcal{X}  =\mathcal{B}$  by means of various tensor-based methods (including BiCG adopted in this paper)  over matricization methods. Please see \cite[Tables 1 and 4]{huang2021numerical} for details of numerical comparison. In particular, if the coefficient tensors $ \mathcal{A}, \mathcal{G}, \mathcal{K}$ of  the ARTE \cref{ARTE}  are of the generalized CPD format  (see \cref{Def: GCPD}), then \cref{alg:1} can save enormous computational costs compared to tensor-matricization methods for  at least two reasons given below. 
\begin{itemize}
    \item If  $ \mathcal{A}, \mathcal{G}, \mathcal{K}$ are of the from of the generalized CPD format, the Einstein products in Steps 1 and 3 of \cref{alg:1} can be performed based on \cref{Pro: 3.12}. For example, the Einstein product $\mathcal{G}*\mathcal{E}_k$ can be computed from the mode products $\sum_{r=1}^R \mathcal{E}_k \times_1 G_1^{(r)} \times_3 G_2^{(r)} \times_5 \cdots \times_{2N-1} G_N^{(r)}$, which only requires $|\mathbf{I}|^2R\sum_{n=1}^N I_n$ multiplications. On the contrary, tensor matricization  methods destroy the inherent tensor structure by transforming  the tensor product $\mathcal{G}*\mathcal{E}_k$ to the matrix product  $\phi(\mathcal{G})\phi(\mathcal{E}_k)$. In this process, firstly, forming the unfolding matrix $\phi(\mathcal{G}) = \sum_{r=1}^R G_N^{(r)}\otimes \cdots \otimes G_2^{(r)} \otimes G_1^{(r)} \in \mathbb{C}^{|\mathbf{I}|\times |\mathbf{I}|}$ from the generalized CPD factor matrices $\{ G_n^{(r)} \}_{r=1,n=1}^{R,N}$  requires  $O(R(N-1)|\mathbf{I}|^2)$ multiplications; secondly, forming the unfolding matrix  $\phi(\mathcal{E}_k)$  requires $O(|\mathbf{I}|)$ multiplications due to the products of tensor indices based on the index mapping `$\mathrm{ivec}$' defined in  equation \cref{eq:ivec}.  Finally, the computation of the matrix product $\phi(\mathcal{G}) \phi(\mathcal{E}_k)\in \mathbb{C}^{|\mathbf{I}|\times |\mathbf{I}|}$ requires $|\mathbf{I}|^3$ multiplications. Adding these together, the computation for $\mathcal{G}*\mathcal{E}_k$  via matricization method requires $O\left(|\mathbf{I}|^3 +R(N-1)|\mathbf{I}|^2\right)$ multiplications. Therefore, the matricization method is much more costly than the tensor-based method.
    
  \item In step 2 of \cref{alg:1}, the involved Einstein product is performed for tensors without specific structure.  When the Newton-VecLTE-BiCG method is used, the Lyapunov tensor equation \cref{17} is converted to the form of multilinear system $\mathcal{A}\ast \mathcal{X}=\mathcal{B}$. It is demonstrated in reference  \cite{huang2021numerical} that the tensor BiCG method (employed in the $\textbf{Newton-VecLTE-BiCG}$ subroutine) is superior to matricization methods on modern computational platforms; see Table 4 in  \cite{huang2021numerical} for detailed comparison.
\end{itemize}

\begin{figure}[htbp!]
		\centering
            \label{fig: reserror}
		\subfloat[residual error]{
			\includegraphics[width=0.45\textwidth]{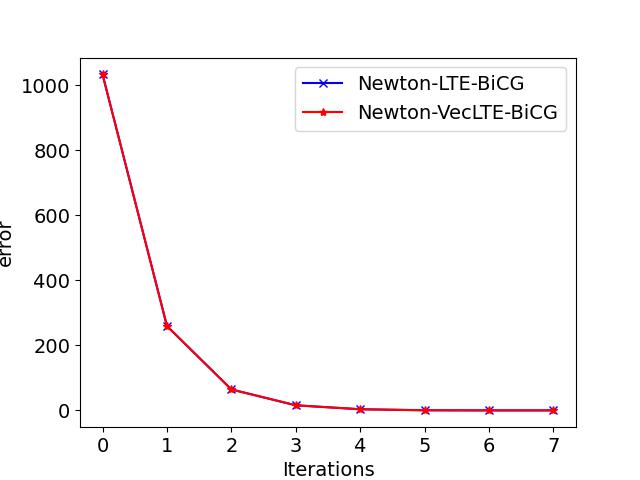}}
		\subfloat[logarithm of the residual error]{
			\includegraphics[width=0.45\textwidth]{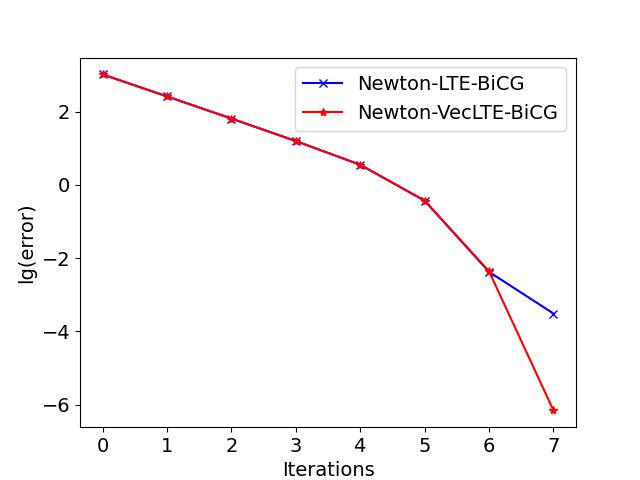}}
		\caption{Residual errors}
\end{figure}

{\bf Example.} Consider the continuous-time MLTI system \cref{MLTI-continuous}, where  the system tensors
$
\mathcal{A} = A_1 \circ A_2 \in \mathbb{C}^{3\times 3 \times 2\times 2}, \mathcal{B} = B_1\circ B_2 \in \mathbb{C}^{3\times 1\times 2\times 1},  \mathcal{C} = C_1\circ C_2 \in \mathbb{C}^{1\times 3\times 1\times 2},
		\mathcal{D} = \mathcal{O}
$
with
\begin{align*}
	\begin{cases}
		A_1 = \begin{bmatrix}
			0 & 1 & 0\\
			0 & 0 & 1\\
			0.2 & 0.5 & 0.8
		\end{bmatrix}, A_2 = \begin{bmatrix}
			0 & 1\\
			0.5 & 0
		\end{bmatrix}, \\
		B_1 = \begin{bmatrix}
			0 &
			0 &
			1
		\end{bmatrix}^{\top}, B_2 = \begin{bmatrix}
			0 &
			1
		\end{bmatrix}^{\top}, 
		C_1 = \begin{bmatrix}
			1 & 0 & 0
		\end{bmatrix}, C_2 = \begin{bmatrix}
			1 & 0
		\end{bmatrix}.
	\end{cases}
\end{align*}
The above matrices are borrowed from \cite[Section 7.1]{chen2021multilinear}, in which it is proved that the controllability and observability tensors are of full unfolding rank. Hence, the tensor pair $(\mathcal{A}, \mathcal{B})$ is stabilizable and $(\mathcal{C}, \mathcal{A})$ is detectable. According to \cref{Thm: ARTEunique}, the ARTE \cref{ARTE} with $\mathcal{G} = \mathcal{B}* \mathcal{B}^{\mathrm{H}} = (B_1B_1^{\mathrm{H}})\circ (B_2B_2^{\mathrm{H}})$ and $\mathcal{K} = \mathcal{C}^{\mathrm{H}}*\mathcal{C} = (C_1^{\mathrm{H}}C_1)\circ (C_2^{\mathrm{H}}C_2)$ has a unique positive semidefinite solution. We perform the tensor products by the built-in function `einsum' in NumPy package \cite{daniel2018opt} with no matricization involved. In the numerical test, we choose the initial tensor $\mathcal{E}_0$ with  
\begin{align*}
	\begin{cases}
		(\mathcal{E}_0)_{::, 11} = \begin{bmatrix}
			10 & 0 & 0 \\
			0 & 4 & 0 \\
			0 & 0 & 13 
		\end{bmatrix},\ 
		(\mathcal {E}_0)_{::,21}  = \begin{bmatrix}
			0 & 0 & 0 \\
			0 & 0 & 0 \\
			1 & 0 & 5\end{bmatrix}, \\
		(\mathcal {E}_0)_{::, 12} = \begin{bmatrix}
			0 & 0 & 1 \\
			0 & 0 & 0\\
			0 & 0 & 5 \end{bmatrix}, \
		(\mathcal {E}_0)_{::, 22}  = \begin{bmatrix}
			7 & 0 & 1\\
			0 & 21 & 5\\
			1 & 5 & 4 \end{bmatrix},
	\end{cases}
\end{align*}
and set the error tolerance $ 10^{-4}$ for the residual norm $\| \mathcal{A}_k^{\mathrm{H}} * \mathcal{E}_{k+1} + \mathcal{E}_{k+1}* \mathcal{A}_k + \mathcal{K}_k \|_F $ of the Lyapunov tensor equation \cref{17} in \cref{alg:1}. \cref{fig: reserror} shows that the Frobenius norm and the logarithm of the Frobenius norm of the error tensor  $f(\mathcal{E}_{k+1})$	of  \cref{eq:f_26Jan} in each Newton iteration. It can be seen that both methods, Newton-LTE-BiCG and Newton-VecLTE-BiCG, yield satisfactory results.
 
	The numerical solution $\mathcal{E} \in \mathbb{C}^{3\times 3 \times 2 \times 2}$ is obtained as given below. Its approximate error is about $\| f(\mathcal{E}) \|_F = 6.9709 \times 10^{-7}$. The smallest U-eigenvalue of $\mathcal{E}$ computed by the higher-order Rayleigh quotient iteration \cite{Brazell2013} is about $0.0063$. Thus, $\mathcal{E}$ is positive definite. The U-eigenvalues of $\mathcal{A}-\mathcal{G}*\mathcal{E}$ are $-1.0165 \pm  0.1846\iota$, $-0.4144 \pm 0.4918\iota$, and $-0.0485 \pm 0.3339\iota$, where $\iota=\sqrt{-1}$. Hence  $\mathcal{A} - \mathcal{G}*\mathcal{E}$ is stable. 
    \small
	\begin{align*}
		\begin{cases}
			\mathcal{E}_{::11}= \begin{bmatrix}
				4.8082 & -0.2001 & 3.9671 \\
				-0.2001 & 1.5958 & -3.3882 \\
				3.9671 & -3.3882 & 18.7381 
			\end{bmatrix}, \ \ 
			\mathcal{E}_{::21}= \begin{bmatrix}
				-0.5391 & -0.0033 & 1.5582 \\
				10.0971 & -4.2223 & 25.4769 \\
				1.1050 & 0.0067 & 5.4633
			\end{bmatrix}, \\
			\mathcal{E}_{::12}= \begin{bmatrix}
				-0.5391 & 10.0971 & 1.1050 \\
				-0.0033 & -4.2223 & 0.0067 \\
				1.5582 & 25.4769 & 5.4633
			\end{bmatrix}, \ \ 
			\mathcal{E}_{::22}= \begin{bmatrix}
				0.9711 & 0.4996 & 0.7895 \\
				0.4996 & 41.7634 & 6.7580 \\
				0.7895 & 6.7580 & 2.9588
			\end{bmatrix}.
		\end{cases}
	\end{align*}
    \normalsize
    
\begin{figure}[htbp!]
		\centering
            \label{fig: conNums}
		\subfloat[The values of $\| \Delta \mathcal{E}\|_F$]{\includegraphics[width=0.44\textwidth]{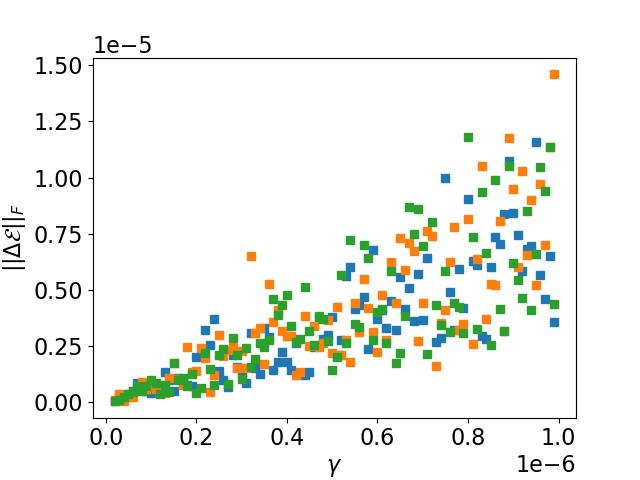}}
		\subfloat[The values of $\frac{\|\Delta \mathcal{E} \|_F}{\Delta_1  \| \mathcal{E} \|_F}$]{\includegraphics[width=0.44\textwidth]{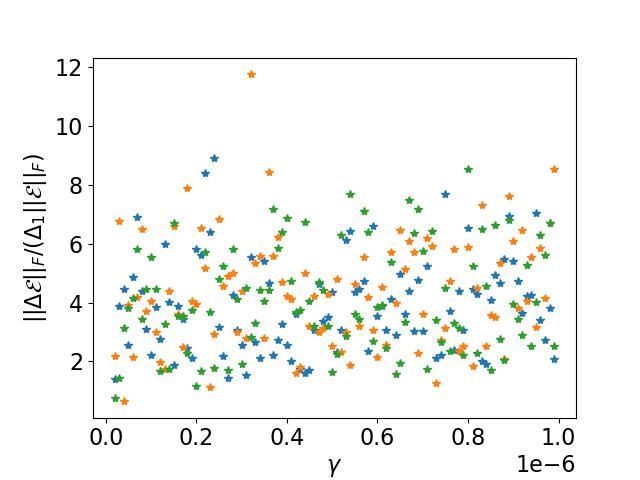}} \\
		\subfloat[The values of $\frac{\|\Delta \mathcal{E} \|_F}{\Delta_2  \| \mathcal{E} \|_F}$]{\includegraphics[width=0.44\textwidth]{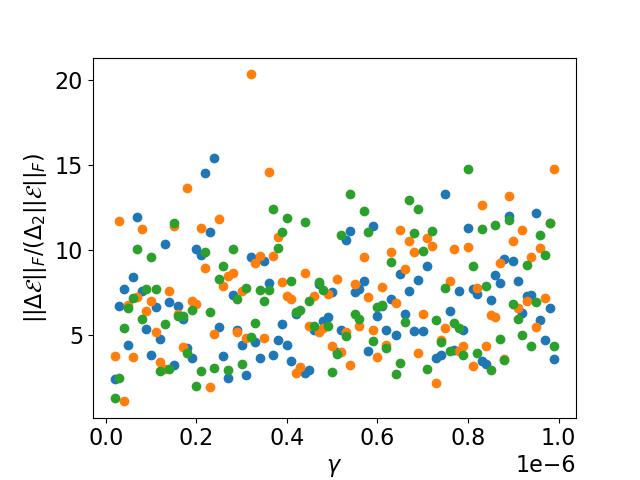}}
		\subfloat[The values of $\frac{\|\Delta \mathcal{E} \|_F}{\Delta_3  \| \mathcal{E} \|_F}$]{\includegraphics[width=0.44\textwidth]{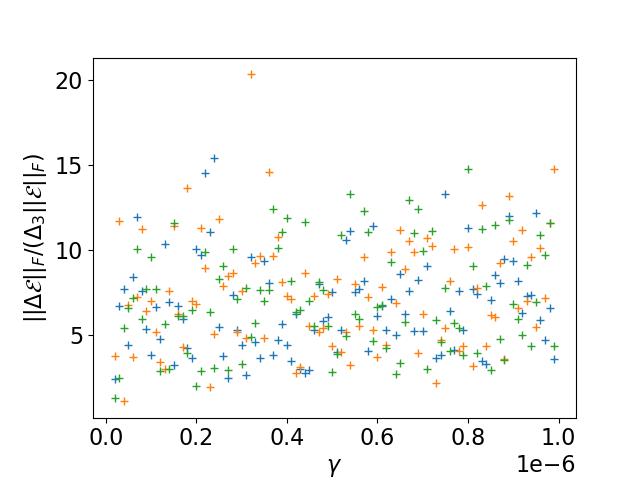}}
		\caption{Perturbation errors.  The randomized experiments for various perturbation factors $\gamma = j\times 10^{-8}\ (j = 1,\cdots, 99)$ are repeated three times, and the numerical results of the three independent experiments are presented in different colors in the four figures.} 
	\end{figure}

	Furthermore, we compute the upper bounds according to \cref{Thm: ConNum} and \cref{Remark 4.2} for the condition numbers; we get  $\kappa_1 \lesssim 55.0299$, $\kappa_2 \lesssim 75.4538$, and $\kappa_3 \lesssim 123.7297$. 
	We conduct the numerical test to check these theoretical upper bounds by randomly setting the perturbations on coefficient tensors. Concretely, $\Delta \mathcal{A} = \gamma \frac{\| \mathcal{A} \|_F }{ \| \Delta \tilde{\mathcal{A}} \|_F}\Delta \tilde{\mathcal{A}}$ where the elements of $\Delta \tilde{\mathcal{A}} $ are independently identically distributed in the standard normal distribution. For the perturbations of $\mathcal{G}$ and $\mathcal{K}$, we add the perturbation only on their nonzero entries, that is, $\Delta \mathcal{G}$ has only one nonzero entry with $\Delta \mathcal{G}_{1111} = \gamma$, and the perturbation $\Delta \mathcal{K}$ has only one nonzero entry with $\Delta \mathcal{K}_{3322} = \gamma$.       We perform randomized experiments independently for $\gamma = j\times 10^{-8}\ (j = 1,\cdots, 99)$ three times. The Frobenius norms $\| \Delta \mathcal{E}\| _F$ are shown in \cref{fig: conNums} (a), and \cref{fig: conNums} (b), (c) and (d) present the magnitudes of $\frac{\|\Delta \mathcal{E} \|_F}{\Delta_1  \| \mathcal{E} \|_F}$, $\frac{\|\Delta \mathcal{E} \|_F}{\Delta_2  \| \mathcal{E} \|_F}$ and $\frac{\|\Delta \mathcal{E} \|_F}{\Delta_3  \| \mathcal{E} \|_F}$, respectively.  It can be seen that they are all smaller than their theoretical  upper bounds given above. \Cref{table} shows the relative errors $\frac{\|\Delta \mathcal{E} \|_F}{\| \mathcal{E} \|_F}$ and the estimates by these condition numbers under three different perturbation tests. Clearly, they all provide valid estimations.
        \begin{table}[htbp!]
		\caption{Comparison of the relative error $\frac{\|\Delta \mathcal{E} \|_F}{\| \mathcal{E} \|_F}$ with the estimates via condition numbers}
		\centering
		\label{table}
		\setlength{\tabcolsep}{5pt}
		\begin{tabular}{c|c|c|c|c}
			\hline
			$\delta$ & 
			$\|\Delta \mathcal{E} \|_F /\| \mathcal{E} \|_F$ & 
			$\kappa_1 \Delta_1$ & $\kappa_2 \Delta_2$ & $\kappa_3 \Delta_3$  \\
			\hline
			$1\times 10^{-8}$ & $4.0547\times 10^{-8}$ & $9.5314 \times 10^{-7} $& $7.5453\times 10^{-7}$ & $1.2373 \times 10^{-6}$ \\
			\hline 
			$2\times 10^{-7}$&$2.2174\times 10^{-6}$ &$1.9062\times 10^{-5}$ & $1.5090\times 10^{-5}$ &$2.4747\times 10^{-5}$ \\
			\hline
			$3\times 10^{-6}$& $3.4478\times 10^{-5}$& $2.8594\times 10^{-4}$&$2.2636\times 10^{-4}$ &$3.7121\times 10^{-4}$ \\
			\hline
		\end{tabular}
		\label{tab1}
	\end{table}

\section{Conclusion}\label{sec:con}
In this paper, we have conducted an in-depth investigation of even-order paired tensors under the Einstein product. For this special class of tensors, we have defined Hamiltonian tensors and symplectic tensors, and we have explored their structural decompositions and properties. We applied these tensor computations to a class of continuous-time multilinear time-invariant control systems and examined a Lyapunov tensor equation. Specifically, we performed a comprehensive study of an algebraic Riccati tensor equation (ARTE) and used it to derive a tensor version of the bounded real lemma and the small gain theorem. For future work, we are interested in developing more efficient numerical algorithms for solving this class of ARTEs. Another important research direction is to extend the current research to the robust control of multilinear control systems. In particular, as a complement to the well-developed gain theory, a new phase theory has recently been proposed for linear control systems (see, for example, \cite{MCQ22, SZQ23}). It would be intriguing to investigate whether phases can be defined for tensors and whether a useful phase theory can be developed for multilinear control systems.

 \section*{Acknowledgments}
 The authors thank the reviewers for their careful reading and constructive comments.

\appendix

\section{Additional tensor basics}

\begin{definition}[Mode product \cite{kolda2009tensor}] \label{def:mode product}
	The mode-$n$ product of a tensor $\mathcal{X} \in \mathbb{C}^{I_1\times I_2 \times \cdots \times I_N}$ and a matrix $A \in \mathbb{C}^{J_n \times I_n}$ is an $I_1 \times \cdots \times I_{n-1} \times J_n \times I_{n+1} \times \cdots \times I_N$ tensor, given by
	\[ \left(\mathcal{X} \times _n A \right)_{i_1\cdots i_{n-1}ji_{n+1}\cdots i_N} = \sum_{i_n=1}^{I_n} \mathcal{X}_{i_1\cdots i_{n-1}i_n i_{n+1}\cdots i_N} A_{ji_n}. \]
\end{definition}

There is a natural connection between even-order paired tensors and matrix outer products. For example, the tensor $\mathcal{A} = A_1\circ A_2 \in \mathbb{C}^{I_1\times J_1 \times I_2\times J_2}$ formed by the outer product of  $A_1 \in \mathbb{C}^{I_1\times J_1}$ and $A_2 \in \mathbb{C}^{I_2\times J_2}$ is a fourth-order paired tensor. In fact, by viewing each paired index as one long index, several classical tensor decompositions such as CANDECOMP/PARAFAC decomposition (CPD) \cite{hitchcock1927expression} and the tensor-train decomposition \cite{oseledets2011tensor} are naturally generalized to even-order paired tensors. Here, we only present the generalized CPD as it is used in \Cref{Section: 3.2} and \Cref{sec. 5}.

\begin{definition}[Generalized CPD \cite{chen2021multilinear}]
	\label{Def: GCPD}
    Let $\mathcal{A} \in \mathbb{C}^{I_1\times J_1\times \cdots \times I_N\times J_N}$ be an even-order paired tensor. Then $\mathcal{A}$ can be decomposed into a sum of a series of matrix outer products, i.e.,
        \begin{equation} \label{GeneralizedCPD}
            \mathcal{A} = \sum_{r=1}^{R} A_1^{(r)}\circ A_2^{(r)}\circ \cdots \circ A_N^{(r)}, 
        \end{equation}
	where $A^{(r)}_n \in \mathbb{C}^{I_n\times J_n}$ for $n=1,\dots,N$ and $r=1,\cdots,R$. The smallest $R$ such that the above equation holds is called the CP rank of $\mathcal{A}$. In the following, the term `rank' of a tensor always refers to the CP rank. 
\end{definition}

Given two even-order paired tensors in the generalized CPD format, their Einstein product can be computed while preserving the structure, thus significantly reducing the computational cost compared to the computation by forming full tensors or tensor matricization.
\begin{proposition}[\cite{chen2021multilinear}]
	\label{Pro: fastGCP}
	Given even-order paired tensors $\mathcal{A}\in \mathbb{C}^{I_1\times J_1 \times \cdots \times I_N\times J_N}$ and $\mathcal{B}\in \mathbb{C}^{J_1\times K_1\times \cdots \times J_N\times K_N}$ in the generalized CPD format of \cref{GeneralizedCPD} with ranks $R$ and $S$, respectively, the Einstein product $\mathcal{A}*\mathcal{B}$ can be computed in an economic fashion as $\mathcal{A}*\mathcal{B} = \sum_{t=1}^{T} F_1^{(t)} \circ F_2^{(t)}\circ \cdots \circ F_N^{(t)}$, 	where $F_n^{(t)} = A_n^{(r)}B_n^{(s)} \in \mathbb{C}^{I_n \times K_n}$ with $t = \mathrm{ivec}((r,s),(R,S))$, and the rank $T=RS$.
\end{proposition}

The following result can be easily established, so its proof is omitted.
\begin{proposition}
	\label{Pro: rank-onetensor}
	For a rank-one even-order paired tensor $\mathcal{A} = A_1 \circ A_2 \circ \cdots \circ A_N \in \mathbb{C}^{I_1\times J_1 \times \cdots \times I_N\times J_N}$,
	\begin{enumerate}
		\item $\mathcal{A}^{\rm H} = A_1^{\rm H} \circ A_2 ^{\rm H} \circ \cdots \circ A_N^{\rm H}$;
		\item if $A_n\  (n=1,\dots, N)$ are upper-triangular matrices, then $\mathcal{A}$ is an upper-triangular even-order paired tensor;
		\item if $A_n \  (n=1,\dots, N)$ are unitary matrices, then $\mathcal{A}$ is a unitary tensor;
		\item  if $A_n \ (n=1,\dots,N)$ are invertible, then $\mathcal{A}$ is also invertible, and $\mathcal{A}^{\rm -1} =  A_1^{-1}\circ A_2^{\rm -1}\circ \cdots \circ A_N^{\rm -1}$.
	\end{enumerate}
\end{proposition}

\section{Proofs} \label{appemd:proofs}
\textbf{The proof of \cref{Pro: TensorKron}.} We prove these properties one by one.
Property 1 can be easily obtained via tensor matricization and \cite[Theorem 3.3]{ragnarsson2012block}.
For Property 2,  firstly, $(\mathcal{A}\otimes \mathcal{B})^{\mathrm{H}} \in \mathbb{C}^{J_1L_1\times I_1K_1 \times \cdots \times J_NL_N \times I_NK_N}$. Then from \cref{TensorKronecker} of tensor Kronecker product, for all $ \mathbf{1} \le \mathbf{i} \le \mathbf{I}, \mathbf{1}\le \mathbf{j}\le \mathbf{J}$, it follows that
\begin{align*}
	\left((\mathcal{A}\otimes \mathcal{B})^{\mathrm{H}}\right)_{[\mathbf{j}\wedge \mathbf{i}]} & = (\mathcal{A}\otimes \mathcal{B})^{\mathrm{H}}( (j_1-1)L_1 + 1: j_1L_1,(i_1 - 1)K_1 +1:i_1K_1, \cdots, 
	\\ & \quad \quad \quad \quad \quad \ \ (j_N - 1)L_N +1 : j_N L_N, (i_N -1)K_N + 1 : i_N K_N)  \\
	& = \overline{\mathcal{A}_{i_1j_1 \cdots i_Nj_N}} \cdot \mathcal{B}^{\mathrm{H}} = (\mathcal{A}^{\mathrm{H}})_{j_1i_1\cdots j_Ni_N}\cdot \mathcal{B}^{\mathrm{H}},
\end{align*}
which shows that $(\mathcal{A}\otimes \mathcal{B})^{\mathrm{H}}$ is the Kronecker product of $\mathcal{A}^{\mathrm{H}}$ and $\mathcal{B}^{\mathrm{H}}$. Next, we prove
Property 3. The tensor $(\mathcal{A}\otimes \mathcal{B})*(\mathcal{C}\otimes \mathcal{D})$ is an $I_1\times R_1\times \cdots \times I_N \times R_N$ block tensor. For $\mathbf{1} \le \mathbf{i} \le \mathbf{I}$, $\mathbf{1} \le \mathbf{r} \le \mathbf{R}$, the $(\mathbf{i}\wedge \mathbf{r})$th block of $(\mathcal{A}\otimes \mathcal{B})*(\mathcal{C}\otimes \mathcal{D})$ is  
\begin{align*}
	& \left( (\mathcal{A}\otimes \mathcal{B})* (\mathcal{C}\otimes \mathcal{D}) \right)_{[\mathbf{i}\wedge \mathbf{r}]} \\
	& \  = \sum_{\mathbf{1}\le \mathbf{j}\le \mathbf{J}} (\mathcal{A}\otimes \mathcal{B})_{[\mathbf{i}\wedge \mathbf{j}]}*(\mathcal{C}\otimes \mathcal{D})_{[\mathbf{j}\wedge \mathbf{r}]}  = \sum_{\mathbf{1}\le \mathbf{j}\le \mathbf{J}} (\mathcal{A}_{i_1j_1\cdots i_N j_N} \cdot \mathcal{B}) * (\mathcal{C}_{j_1r_1\cdots j_N r_N} \cdot \mathcal{D}) \\
	& \ =  \mathcal{B} *\mathcal{D}\sum_{\mathbf{1}\le \mathbf{j}\le \mathbf{J}}\mathcal{A}_{i_1j_1\cdots i_N j_N} \mathcal{C}_{j_1r_1\cdots j_N r_N}= (\mathcal{A}*\mathcal{C})_{i_1r_1\cdots i_Nr_N}\cdot (\mathcal{B} *\mathcal{D})  \\
	& \ =\left((\mathcal{A}* \mathcal{C}) \otimes (\mathcal{B}* \mathcal{D})\right)_{[\mathbf{i}\wedge \mathbf{r}]}.
\end{align*}
Thus, $(\mathcal{A}\otimes \mathcal{B}) * (\mathcal{C}\otimes \mathcal{D}) = (\mathcal{A}*\mathcal{C})\otimes (\mathcal{B}* \mathcal{D})$.  Property 4 follows from Property 3. $\Box$

\textbf{The proof of \cref{Pro: KronVec}.} 	We know that both $\mathrm{Vec}(\mathcal{U}* \mathcal{X}* \mathcal{W})$ and $(\mathcal{W}^{\top} \otimes \mathcal{U})* \mathrm{Vec}(\mathcal{X}) \in \mathbb{C}^{L_1I_1\times \cdots \times L_NI_N}$ are $L_1\times \cdots \times L_N$ block tensors, each subblock is of size $I_1 \times \cdots \times I_N$. Hence, it suffices to compare all subblock tensors.  For $\mathbf{l}=(l_1,\dots, l_N)$, the $\mathbf{l}$-th subblock $\mathrm{Vec}(\mathcal{U}* \mathcal{X}* \mathcal{W})_{[\mathbf{l}]} = \left(\mathcal{U}* \mathcal{X}* \mathcal{W}\right)_{:l_1, \cdots ,:l_N}  = \mathcal{U}* \mathcal{X}* \mathcal{W}_{:l_1, \cdots ,:l_N}$, and $\left((\mathcal{W}^{\top} \otimes \mathcal{U})* \mathrm{Vec}(\mathcal{X})\right)_{[\mathbf{l}]} =  \sum_{\mathbf{1}\le \mathbf{k}\le \mathbf{K}} (\mathcal{W}^{\top}\otimes \mathcal{U})_{[\mathbf{l}\wedge \mathbf{k}]} * \mathrm{Vec}(\mathcal{X})_{[\mathbf{k}]}  = \sum_{\mathbf{1}\le \mathbf{k}\le \mathbf{K}} \mathcal{W}_{k_1l_1\cdots k_Nl_N}\mathcal{U}*\mathcal{X}_{:k_1,\cdots, : k_N} = \mathcal{U}* \mathcal{X}* \mathcal{W}_{:l_1,\cdots, :l_N}$. Thereby, the result holds.

\section{Fast operations for structured tensors} \label{appendix C}

The following result can be easily checked.
\begin{proposition} 
	\label{Pro: 3.12}
	Given $\mathcal{E}\in \mathbb{C}^{I_1 \times J_1 \times  \cdots \times I_N\times J_N}$ and a tensor $\mathcal{A} = \sum_{r=1}^R A_1^{(r)} \circ \cdots \circ A_N^{(r)}$ in generalized CPD format, we have
	\begin{enumerate}
		\item $\mathcal{E}*\mathcal{A} = \sum_{r=1}^R \mathcal{E} \times_2 A_1^{(r)\top} \times_4 A_2^{(r)\top} \times_6 \cdots \times_{2N} A_N^{(r)\top}$ for $A_n^{(r)} \in \mathbb{C}^{J_n \times K_n}$,
		\item $\mathcal{A}*\mathcal{E} = \sum_{r=1}^{R} \mathcal{E} \times_1 A_1^{(r)} \times_3 A_2^{(r)} \times_5 \cdots \times_{2N-1} A_N^{(r)}$ for $A_n^{(r)} \in \mathbb{C}^{K_n \times I_n}$.
	\end{enumerate} 
\end{proposition} 
\begin{remark}
    It requires $\left(|\mathbf{I}| + R(N -1)\right)|\mathbf{J}| |\mathbf{K}|$ multiplications for the above Einstein product $\mathcal{E}*\mathcal{A}$ based on full tensor or tensor matricization, but only needs $R|\mathbf{I}| |\mathbf{K}| (\sum_{n=1}^N J_n)$ multiplications according to the equations in \cref{Pro: 3.12}, thus saving a lot of computational and storage costs. 
\end{remark}

\begin{remark}\label{rem:lyap}
    From \cref{Pro: 3.12,Pro: rank-onetensor}, the ARTE \cref{ARTE} reduces to the Lyapunov tensor equation  
    \[
		\mathcal{E}\times_1 A + \mathcal{E} \times_2 A + \cdots + \mathcal{E} \times_{2N} A + \mathcal{K} = \mathcal{O}
  \] with $A \in \mathbb{R}^{m\times m}$,
	provided that 
	$\mathcal{G}=\mathcal{O}$ and $\mathcal{A}$ have the generalized CPD format
	\begin{equation} \label{eq:tensorA}
		\mathcal{A} = A^{\top}\circ I_m\circ\cdots \circ I_m + I_m\circ A^{\top} \circ \cdots \circ I_m + \cdots + I_m\circ \cdots \circ I_m\circ A^{\top} . 
	\end{equation}
\end{remark}

\begin{proposition}
	\label{Pro: KronSum}
	Let $\mathcal{A}\in \mathbb{C}^{m\times m \times \cdots \times m\times m}$ be the $2N$th-order paired tensor as given in \cref{eq:tensorA}, and the eigenvalues of $A$ are $\{ \lambda_1, \lambda_2, \dots, \lambda_m\}$. Then the U-eigenvalues of $\mathcal{A}$ are $\Lambda(\mathcal{A}) = \{ \lambda_{i_1}+\cdots +\lambda_{i_N} \big| i_n=1,\dots,m \ \mathrm{for} \ \mathrm{all}\ n=1,\dots,N \}$. 
\end{proposition}
\begin{proof}
	Let the Schur decomposition of the matrix  $A^{\top}$ be  $A^{\top} = URU^{\rm H}$,
	where $U\in \mathbb{C}^{m\times m}$ is a unitary matrix, and $R\in \mathbb{C}^{m\times m}$ is an upper-triangular matrix whose diagonal entries $R_{kk} = \lambda_k, k=1,\dots,m$. For all $n=1,\dots,N$ and the rank-one tensor $ \mathcal{A}_n = \underbrace{ I_m\circ \cdots \circ I_m }_{n-1\ \mathrm{terms}}\circ A^{\top} \circ \underbrace{I_m\circ \cdots\circ I_m}_{ N-n\ \mathrm{terms}}$, it follows from \cref{Pro: fastGCP,Pro: rank-onetensor} that $\mathcal{A}_n  = \left( UI_mU^{\rm H} \right)\circ \cdots \circ \left( UI_mU^{\rm H} \right) \circ \left( URU^{\rm H} \right)\circ \left( UI_mU^{\rm H} \right) \circ \cdots \circ \left( UI_mU^{\rm H} \right)  = \mathcal{U}*\mathcal{R}_n*\mathcal{U}^{\rm H}$, where $\mathcal{U} = U\circ U \circ \cdots \circ U \in \mathbb{C}^{m\times m \times \cdots \times m \times m}$ is a unitary tensor, and  $\mathcal{R}_n =  I_m\circ \cdots \circ I_m \circ R \circ I_m\circ \cdots\circ I_m$ is an upper-triangular $2N$th-order paired tensor. Hence,  $\mathcal{A} = \mathcal{U}*\mathcal{R}*\mathcal{U}^{\rm H}$ with $\mathcal{R} = \left( \mathcal{R}_1+\cdots +\mathcal{R}_N \right) $ is the tensor Schur decomposition of $\mathcal{A}$. According to \cref{lem:schur}, the diagonal entries of $\mathcal{R}$ are the U-eigenvalues of $\mathcal{A}$. From the definition of the  outer product, the diagonal elements, i.e., the U-eigenvalues of $\mathcal{A}$ are $\mathcal{R}_{i_1i_1\cdots i_Ni_N} = (\mathcal{R}_1)_{i_1i_1\cdots i_Ni_N} + \cdots + (\mathcal{R}_N)_{i_1i_1\cdots i_Ni_N} = R_{i_1i_1} + R_{i_2i_2} + \cdots + R_{i_Ni_N} = \lambda_{i_1} + \lambda_{i_2} + \cdots + \lambda_{i_N}$ for all $i_n = 1,\dots, m$ and $n=1,\dots,N$.
\end{proof}

\begin{remark}
	According to \cref{Pro: KronSum}, the tensor $\mathcal{A}$ in  \cref{eq:tensorA}  is  stable  if and only if  
	the matrix $A$ is  Hurwitz stable.
\end{remark}

For rank-one even-order paired tensors, the following results can be easily established; therefore, their proofs are omitted.
\begin{proposition}
    \label{Pro: rankoneSchur}
    Let $\mathcal{A} = A_1 \circ A_2 \circ \cdots \circ A_N$ be a rank-one even-order paired tensor, where $A_k \in \mathbb{C}^{I_k\times I_k}$. Given the Schur decompositions $A_n = U_nR_nU_n^{\rm H}$ for all $n= 1,\dots,N$, the tensor Schur decomposition of $\mathcal{A}$ is
	$ \mathcal{A} = \mathcal{U}*\mathcal{R}*\mathcal{U}^{\rm H}$, 
	where $\mathcal{U} = U_1\circ U_2 \circ \cdots \circ U_N$ is unitary, and $\mathcal{R} = R_1\circ R_2\circ \cdots \circ R_N$ is upper-triangular. Therefore, the U-eigenvalues of $\mathcal{A}$ are the products of the eigenvalues of $A_n\ (n=1,\dots,N)$, i.e., $\Lambda(\mathcal{A}) = \{ \prod_{n=1}^{N} \lambda_n \big| \ \lambda_n \in \Lambda(A_n) \}$.
\end{proposition}

\begin{proposition}
    \label{Pro: rankonenorm}
	For the rank-one even-order paired tensor $\mathcal{A} = A_1\circ \cdots \circ A_N$, $
		\| \mathcal{A} \|_{\star}  =  {\textstyle \prod_{n=1}^{N}} \| A_n \|_{\star}$, 
	where $\star $ denotes either the Frobenius norm or the spectral norm. 
\end{proposition}

\begin{proposition} \label{prop:kron}
	$
		\| \mathcal{A} \otimes \mathcal{B} \|_{\star}  =  \| \mathcal{A} \|_{\star} \| \mathcal{B} \|_{\star}$, 
	where $\star $ denotes either the Frobenius norm or the spectral norm.
\end{proposition}

\bibliographystyle{siamplain}
\bibliography{ref}

\end{document}